\documentclass[10pt,draft,twoside,leqno,a4paper]{article}     
\usepackage{amsmath,amsthm,amsfonts,amssymb,verbatim}     
\newtheorem{Theorem}{Theorem}[section]    
\newtheorem{Lemma}{Lemma}[section]    
\newtheorem{Corollary}{Corollary}[section]    
\newtheorem{Proposition}{Proposition}[section]    
\newtheorem{Remark}{Remark}[section]   
\newtheorem{Question}{Question}[section]   
\newcommand{\R}{\mathbb R}     
\newcommand{\eps}{\varepsilon}     
\newcommand{\de}{\partial}     
\newcommand{\intRn}{\int_{\mathbb R^n}}     
\newcommand{\dx}{\,dx}     
\newcommand{\dy}{\,dy}     
\newcommand{\weak}{\rightharpoondown}     
\newcommand{\bra}{\langle}     
\newcommand{\ket}{\rangle}   
\newcommand{\br}{\bar r}    
\newcommand{\Ia}{I_\alpha}    
    
\newcommand{\intM}{\int_M}     
\newcommand{\Delg}{\Delta_g}     
\newcommand{\dvg}{\,dv_g}     
\newcommand{\nabg}{\nabla_g}     
\newcommand{\dsg}{\,ds_g}     
\newcommand{\Rg}{R_g}     
\newcommand{\Aq}{\mathcal A_q}     
\newcommand{\uq}{u_q}

\newcommand{\Yg}{Y_g}     
\newcommand{\Ga}{G_\alpha}     
     
\newcommand{\gh}{\widehat g_\alpha}     
\newcommand{\dvgh}{\,dv_{\widehat g_\alpha}}     
\newcommand{\nabgh}{\nabla_{\widehat g_\alpha}}     
\newcommand{\Delgh}{\Delta_{\widehat g_\alpha}}     
\newcommand{\zea}{\zeta_\alpha}     
\newcommand{\Vatilde}{\widetilde V_\alpha}     
\newcommand{\delxa}{\delta_{x_\alpha}}     
\newcommand{\distg}{{\rm dist}_g}     
     
\newcommand{\volg}{{\rm vol}_g}     
\newcommand{\stpr}{{2^*}'}      
\newcommand{\xa}{x_\alpha}     
\newcommand{\Oa}{\Omega_\alpha}
\newcommand{\tOa}{\widetilde \Omega_\alpha}     
\newcommand{\ga}{g_\alpha}     
\newcommand{\inta}{\int_{\Omega_\alpha}}
\newcommand{\intta}{\int_{\widetilde \Omega_\alpha}}     
     
\newcommand{\dvga}{\,dv_{g_\alpha}}     
\newcommand{\Delga}{\Delta_{g_\alpha}}     
\newcommand{\nabga}{\nabla_{g_\alpha}}     
     
\newcommand{\ua}{u_\alpha}     
\newcommand{\mua}{\mu_\alpha}     
\newcommand{\epsa}{\varepsilon_\alpha}     
\newcommand{\va}{v_\alpha}     
\newcommand{\wa}{w_\alpha}     
\newcommand{\ha}{h_\alpha}     
\newcommand{\xia}{\xi_\alpha}     
\newcommand{\sia}{\sigma_\alpha}     
\newcommand{\la}{\ell_\alpha}     
\newcommand{\ta}{t_\alpha}     
\newcommand{\chia}{\chi_\alpha}     
     
\newcommand{\phia}{\varphi_\alpha}     
\newcommand{\ya}{\tilde x_\alpha}     
\newcommand{\lda}{\lambda_\alpha}     
\newcommand{\Wa}{W_\alpha}

\newcommand{\Ba}{B_\alpha}
\newcommand{\tBa}{B_\alpha}

\newcommand{\barl}{\bar\lambda}     
     
\newcommand{\da}{\delta_\alpha}

\newcommand{\beps}{\bar\varepsilon}     
\newcommand{\Ha}{H_\alpha}     
   
\newcommand{\ra}{r_\alpha}     
\newcommand{\hg}{\widehat g}

\newcommand{\nabhg}{\nabla_{\widehat g}} 
\newcommand{\dvhg}{\,dv_{\widehat g}} 
\newcommand{\dshg}{\,ds_{\widehat g}} 
 
\newcommand{\Delhg}{\Delta_{\widehat g}} 
\newcommand{\intBa}{\int_{B_\alpha}}
\newcommand{\inttBa}{\int_{B_\alpha}} 
\newcommand{\Tmua}{T_{\mu_\alpha}} 
\newcommand{\ka}{k_\alpha} 
\newcommand{\Thetaa}{\Theta_\alpha} 
\newcommand{\fa}{f_\alpha} 
\newcommand{\Qa}{Q_\alpha}    
\begin{document}    
\title{A sharp Sobolev inequality\\on Riemannian manifolds}     
\author{YanYan Li\thanks     
{Partially supported by National Science Foundation  
Grant DMS-0100819 and a Rutgers University Research Council Grant.}\\  
\small{Department of Mathematics}\\  
\small{Rutgers University}\\           
\small{110 Frelinghuysen Rd.}\\ 
\small{Piscataway, NJ 08854-8019, U.S.A.}\\   
\small{yyli@math.rutgers.edu}\\  
\small{http://www.math.rutgers.edu/\~{}yyli}\\       
\and \\     
Tonia Ricciardi\thanks     
{Partially supported by CNR Fellowship 203.01.69 (19/01/98)     
and by PRIN~2000 ``Variational Methods and Nonlinear Differential Equations".}\\   
\small{Dipartimento di Matematica e Applicazioni}\\             
\small{Universit\`a di Napoli Federico II}\\     
\small{Via Cintia}\\ 
\small{80126 Naples, Italy}\\ 
\small{tonia.ricciardi@unina.it}\\ 
\small{http://cds.unina.it/\~{}tonricci}\\         
}   
\date{January 18, 2002}  
\maketitle     
\begin{abstract}     
Let $(M,g)$ be a smooth compact Riemannian manifold      
without boundary of dimension $n\ge6$.      
We prove that      
\begin{align*}      
\|u\|_{L^{2^*}(M,g)}^2    
\le K^2\intM\big\{|\nabg u|^2+c(n)R_gu^2\big\}\dvg      
+A\|u\|_{L^{2n/(n+2)}(M,g)}^2,      
\end{align*}      
for all $u\in H^1(M)$, where    
$2^*=2n/(n-2)$, $c(n)=(n-2)/[4(n-1)]$,      
$R_g$ is the scalar curvature,    
$K^{-1}=\inf\|\nabla u\|_{L^2(\R^n)}\|u\|_{L^{2n/(n-2)}(\R^n)}^{-1}$     
and $A>0$ is a constant     
depending on $(M,g)$ only.      
The inequality is {\em sharp} in the sense that on any $(M,g)$,      
$K$ can not be replaced by any smaller number and $R_g$      
can not be replaced by any continuous function which is   
smaller than $R_g$ at some point.      
If $(M,g)$ is not locally conformally flat,    
the exponent $2n/(n+2)$ can not be replaced by any smaller number.     
If $(M,g)$ is locally conformally flat,  
a stronger inequality, with $2n/(n+2)$ replaced by $1$,  
holds in all dimensions $n\ge 3$.    
\end{abstract}   
\noindent  
{\sc key words:} sharp Sobolev inequality, critical exponent,   
Yamabe problem  
\setcounter{section}{-1}     
\section{Introduction}     
\label{secintro}  
Considerable work has been devoted to   
the analysis of sharp Sobolev-type inequalities,  
very often in connection with concrete problems  
from geometry and physics. See, e.g.,   
Trudinger~\cite{Tr}, Moser~\cite{M},  
Aubin~\cite{Au1,Au2},  Talenti~\cite{Ta}, 
Brezis and Nirenberg~\cite{BN}, Lieb~\cite{Lieb2},
Carleson and Chang~\cite{CC},   
Struwe~\cite{struwe}, Escobar~\cite{E3}, and Beckner~\cite{Beckner}.  

In order to fix notation, we recall 
a classical result.
For $n\ge 3$ and $2^*=2n/(n-2)$,
it was shown by Aubin~\cite{Au1} and Talenti~\cite{Ta} that
\begin{equation}
\label{Rnsobolev}
K^{-1}=\inf \left\{
\frac{\|\nabla u\| _{L^2(\R^n)}}
{\|u\|_{L^{2^*}(\R^n)}}
\ :\
u\in  L^{2^*}(\R^n)\setminus\{0\},
\ |\nabla u|\in  L^2(\R^n)\right\},
\end{equation}
where $K^2=4/[n(n-2)\sigma_n^{2/n}]$ and where $\sigma_n$ is the volume
of the standard $n$-sphere.
They also showed that the infimum is attained and, modulo
non-zero constant multiples, the
set of minimizers is given by
\[
\{U_{y,\lambda}\ ;\
y\in\R^n,\lambda>0\}
\]
where
\begin{align*}
&U_{y,\lambda}(x)=\lambda^{(n-2)/2}U(\lambda(x-y))\\
&U(x)=U_{0,1}(x)=\big(\frac{1}{1+\barl^2|x|^2}\big)^{\frac{n-2}{2}}
\end{align*}
and $\barl^2=[n(n-2)]^{-1}K^{-2}$.
The function $U$ is characterized as the unique solution of
the equation
\begin{equation}
\label{bubbleeq}
-\Delta U=K^{-2}U^{2^*-1}\qquad\text{in}\;\R^n
\end{equation}
satisfying
\begin{align*}
&U\in D^{1,2}(\R^n),\quad0<U\le1,\\
\nonumber
&U(0)=1, \quad\int_{\R^n}U^{2^*}\,dx=1.
\end{align*}

 A conjecture was made  by Aubin~\cite{Au1}:
On any smooth compact Riemannian manifold $(M,g)$ of
dimension $n\ge 3$,  there
exists a constant $A>0$ depending  only on $(M,g)$, such that
\begin{equation}
\label{hebeyvaugon}
\|u\|_{L^{2^*}(M,g)}^2\le K^2\|\nabg u\|_{L^2(M,g)}^2+A\|u\|_{L^2(M,g)}^2,
\qquad \forall\ u\in H^1(M).
\end{equation}
The conjecture was proved  in \cite{Au1}
 for manifolds of constant
sectional curvature. He also proved a weaker version of
\eqref{hebeyvaugon}, where for any $\eps>0$,
$K$ is replaced by $K+\eps$ and
where $A$ is allowed to depend on $\eps$.

Various related questions  in bounded domains $\Omega$
of $\R^n$ have been extensively studied.
In particular, the following result
was proved by Brezis and Nirenberg~\cite{BN}:
For 
$n=3$, there exists   
a constant $\lambda^*>0$ such that   
\[   
\|\nabla u\|_{L^2(\Omega)}^2\ge K^{-2}\|u\|_{L^6(\Omega)}^2  
+\lambda^*\|u\|_{L^2(\Omega)}^2,\qquad   
\forall\ u\in H^1_0(\Omega),   
\]  
where $\lambda^*$ depends on $\Omega$;     
when $\Omega$ is a ball,    
$\lambda^*$ can be taken as $  \frac 14 \pi^2  
(3|\Omega|/(4\pi))^{-2/3}$ which is sharp.   
They also showed that   
for $n\ge 4$ and for all $q<n/(n-2)$,   
\[   
\|\nabla u\|_{L^2(\Omega)}^2\ge K^{-2}\|u\|_{L^{2^*}(\Omega)}^2+\lambda_q   
\|u\|_{L^q(\Omega)}^2,\qquad   
\forall\ u\in H^1_0(\Omega),   
\]   
where $\lambda_q>0$ depends on $\Omega$ and $q$.   
On the other hand they pointed out that, on   
any $\Omega$,     
such an inequality can not hold with $q=n/(n-2)$.   
Further results were obtained by Brezis and Lieb~\cite{BL},   
and closely related ones  by Adimurthi and Yadava~\cite{AY}.    
Results of similar nature, concerning the  
Hardy-Littlewood inequality for   
functions with support in a ball of $\R^n$, were   
obtained and used by Daubechies and Lieb~\cite{DL}.   
We refer to Brezis and Marcus~\cite{BM},   
Brezis, Marcus and Shafrir~\cite{BMS},   
and Shafrir \cite{Sh} for more recent related works on the sharp  
Hardy-Littlewood inequality.

The conjecture (\ref{hebeyvaugon}) was  proved
by Hebey and Vaugon \cite{HV}.
Results of similar nature for manifolds with boundary   
were established by Li and Zhu~\cite{LZ,LZ3}, with improvements
 given by Zhu~\cite{Zhu1,Zhu2}.   
A $W^{1,p}$ version of \eqref{hebeyvaugon} with    
$p\ne 2$, also conjectured by Aubin~\cite{Au1}, was   
proved through the work of Aubin and Li \cite{AL},   
and Druet~\cite{D1,D2}.   
It should be mentioned that in two dimensions,  
the corresponding inequality   
discovered by Trudinger~\cite{Tr} has also  
been widely investigated and applied in its sharp form,  
due to Moser~\cite{M}.  
A sharp Moser-Trudinger inequality on Riemannian 2-manifolds was established   
by Fontana~\cite{F}, and used by     
Gillet and Soul\'e~\cite{GS}. Alternative   
proofs of Fontana's result and connections to the analysis of   
vortices in the Chern-Simon-Higgs gauge theory were given by   
Ding, Jost, Li and Wang~\cite{DJLW} and   
by Nolasco and Tarantello~\cite{NT1,NT2}.\par     
\paragraph{Statement of the main results}  
Our main result in this paper is the following 
sharp Sobolev inequality on Riemannian manifolds 
of dimension $n\ge 6$:
\begin{Theorem} [Main Result]    
\label{thmmain}    
Let $(M,g)$ be a smooth  compact      
Riemannian manifold without boundary  
of dimension    
$n\ge6$.                                  
There exists a constant $A>0$, depending on $(M,g)$ only,    
such that for all $u\in H^1(M)$ there holds:    
\begin{align}    
\label{optsobolev}    
\|u\|_{L^{2^*}(M,g)}^2    
\le K^2\intM\big\{|\nabla_g u|^2+c(n)R_g u^2\big\}\dvg    
+A\|u\|_{L^{\br}(M,g)}^2,    
\end{align}    
where $2^*$ and $K$ are defined above,     
$c(n)=(n-2)/[4(n-1)]$, $\br=2n/(n+2)={2^*}'$,    
$\Rg$ is the scalar curvature of $g$.    
\end{Theorem} 
We point out  
that our proof of Theorem \ref{thmmain}   
does not make any use of inequality \eqref{hebeyvaugon}, 
which on the other hand is an easy consequence.
\begin{Remark} [Sharpness]  
\label{remsharp}   
Theorem \ref{thmmain} is {\em sharp},     
in the sense that    
one can neither replace $K$ by     
any smaller number, nor replace    
$R_g$ by any $R_g+f$ with $f\in C^0$ negative somewhere.    
Moreover, if $(M,g)$ is not locally conformally flat,    
one cannot replace $\br$ by any smaller number.    
\end{Remark}   
The case of locally conformally flat manifolds is  
completely described by our next result:   
\begin{Theorem}   
\label{thmflat}   
Let $(M,g)$ be a smooth compact locally conformally flat   
Riemannian manifold without boundary of  
dimension $n\ge 3$.   
There exists a constant $A>0$, depending on $(M,g)$ only,   
such that for all $u\in H^1(M)$ there holds:   
\begin{align}   
\label{opt}   
\|u\|_{L^{2^*}(M,g)}^2   
\le K^2\intM\big\{|\nabla_g u|^2+c(n)R_g u^2\big\}\dvg   
+A\|u\|_{L^1(M,g)}^2.   
\end{align}   
\end{Theorem}   
In view of the work of Schoen \cite{Schoen},    
we expect   
a positive answer to    
\begin{Question}  
\label{quesglobal}   
For locally conformally flat   
manifolds and for manifolds of dimension $3\le n\le 5$,    
are there some Sobolev type inequalities   
involving {\em global} geometric quantities?   
\end{Question}   
For manifolds with positive total scalar curvature   
$\intM R_g\dvg$, a natural    
global geometric quantity is the ``mass'', which 
corresponds to   
the leading term of the regular part of the   
Green's function for  
the conformal Laplacian, see \cite{Schoen, LP}.\par  
The {\em sharpness} of Theorem \ref{thmmain} as stated   
in Remark \ref{remsharp} can  be deduced from the following   
expansions   
due to Aubin \cite{Au2} (see also \cite{Au,LP}).    
Let $P\in M$, $\lambda>0$ and let $h$ be a Riemannian metric on $M$; 
denote by $\xi_{P,\lambda}^h$ 
the ``$h$-bubble" defined for $x\in M$ by 
\[ 
\xi_{P,\lambda}^h(x)=\left(\frac{\lambda} 
{1+(\lambda\barl)^2\mathrm{dist}_h^2(x,P)}\right)^{\frac{n-2}{2}}. 
\]  
Let $\eta$ be a smooth  cutoff function supported near $P$,   
and set 
$\widetilde\xi_{P,\lambda}^h=\eta\,\xi_{P,\lambda}^h$.  
Then, as $\lambda\to\infty$,    
\begin{align}   
\label{aubinexpansion}    
Y_h(\widetilde \xi_{P,\lambda}^h)=    
\begin{cases}    
K^{-2}-\gamma_n|W_h(P)|^2\lambda^{-4}+\circ(\lambda^{-4}),    
&\text{if}\;n\ge 7\\    
K^{-2}-\gamma_n|W_h(P)|^2\lambda^{-4}\log\lambda+\circ(\lambda^{-4}\log\lambda),    
&\text{if}\;n=6,    
\end{cases}    
\end{align}    
where $\gamma_n>0$ is a dimensional   
constant,  $W_h(P)$ is the Weyl tensor of $h$ at $P$  and   
$Y_h$ denotes the Yamabe functional:  
\begin{align}   
\label{Yamabe}   
Y_h(u)=\frac{\intM\{|\nabla_h u|^2+c(n)R_h u^2\}\,dv_h}   
{\big(\intM|u|^{2^*}\,dv_h\big)^{2/2^*}},   
\qquad u\in H^1(M)\setminus\{0\}.  
\end{align}   
To see the sharpness of Theorem \ref{thmmain}, 
we note that if $K$ is replaced by any smaller number,   
then \eqref{optsobolev} is violated   
by $u=\widetilde \xi_{P,\lambda}^g$ for large $\lambda$   
(fixing any $P\in M$); if   
$R_g$ is replaced by $R_g+f$ with $f(\bar P)<0$ for some $\bar P\in M$, then   
\eqref{optsobolev} is violated   
by $u=\widetilde \xi_{\bar P,\lambda}^g$ for large $\lambda$;   
if $\br$ is replaced by some $1\le s<\br$, then   
we have   
$\|\widetilde \xi_{P,\lambda}^g\|_{L^s(M)}  
=\circ(\|\widetilde \xi_{P,\lambda}\|_{L^{\br}(M)})$,   
and thus inequality \eqref{optsobolev}  and   
\eqref{aubinexpansion}--\eqref{Yamabe} imply that $|W_g(P)|=0$ for all   
$P\in M$, i.e., $(M,g)$ is locally conformally flat.\par  
In view of  \eqref{aubinexpansion} and our results   
we  expect a positive answer to    
\begin{Question}Are there some refined versions of   
\eqref{optsobolev} involving the Weyl tensor?   
\end{Question}   
\paragraph{Outline of the proofs}  
We first sketch   
the proof of Theorem \ref{thmflat},  
which is simple, and    
relies on a ``local to global" argument,  given
in the Appendix.  
By a local to global argument, we mean that   
we first establish the inequality  for all   
functions $u\in H^1(M)$ supported in a ball of  
fixed diameter $\eps>0$, and then    
we extend the inequality to arbitrary $u\in H^1(M)$.   
Another ingredient is the following well-known transformation property 
of the conformal Laplacian,  
see, e.g.,  \cite{SY}:  
\begin{equation}  
\label{conftransf}  
-\Delta_{\widehat h}u+c(n)R_{\widehat h}u 
=\varphi^{1-2^*}\{-\Delta_h(u\varphi)+c(n)R_h(u\varphi)\}, 
\end{equation}  
for all $u\in H^1(M)$, where $\widehat h=\varphi^{4/(n-2)}h$,   
$\varphi\in C^\infty(M)$, $\varphi>0$.  
\begin{proof} [Proof of Theorem \ref{thmflat}]   
Since $(M,g)$ is locally conformally flat, for   
some $\varepsilon>0$ independent of $P\in M$,  
we have $(B_\eps(P),g)\cong(B,\varphi^{4/(n-2)}\mathcal E)$, 
for some $\varphi>0$ (under control), with $\mathcal E$ the Euclidean metric. 
Since   
\[  
(\int_B|u|^{2^*}\dx)^{2/2^*}\le K^2\int_B|\nabla u|^2\dx,\qquad   
\forall u\in H^1_0(B),   
\]   
we have by \eqref{conftransf},    
\[   
(\int_{B_\eps(P)}|u|^{2^*}\dvg)^{2/2^*}\le K^2   
\int_{B_\eps(P)}\{|\nabg u|^2   
+c(n)\Rg u^2\}\dvg, \ \forall u\in H^1_0( B_\eps(P)).   
\]   
Now Theorem \ref{thmflat} follows from the above    
and from Lemma \ref{lemlocaltoglobal} in the Appendix.   
\end{proof}   
The ``local to global" approach has been systematically used by   
Aubin~\cite{Au1}, Hebey and Vaugon~\cite{HV},   
Aubin and Li~\cite{AL},    
Druet, Hebey and Vaugon~\cite{DHV2},    
and others.   
In \cite{LZ,LZ3}, Li and Zhu  
introduced a  global approach  
by attacking the problem directly   
on the whole manifold.  Such an approach   
should be useful in obtaining a positive answer   
to Question \ref{quesglobal}, since the inequality would involve   
global quantities and therefore could not be  
obtained   
by a local to global approach.\par   
We shall now provide a brief sketch of the proof of Theorem \ref{thmmain},  
which will occupy the main part of this paper.  
For simplicity of exposition,   
we shall restrict ourselves in the present sketch to the case $n\ge7$.  
We argue by contradiction,  
and we take a  global approach. Namely,   
for all $\alpha>0$ we define:  
\[  
\Ia(u)=  
\frac{\intM\{|\nabg u|^2+c(n)\Rg u^2\}\dvg+\alpha\|u\|_{L^{\br}(M,g)}^2}  
{\|u\|_{L^{2^*}(M,g)}^2},  
\qquad u\in H^1(M)\setminus\{0\}.  
\]  
Negating \eqref{optsobolev}, we assume that   
\begin{equation}  
\label{introcontradfirst}  
\inf_{H^1(M)\setminus\{0\}}\Ia<K^{-2},  
\quad\forall \alpha>0.  
\end{equation}  
It is straightforward to check that inequality \eqref{optsobolev}  
holds for the family   
$\{t\widetilde\xi_{P,\lambda}^g\}$ defined above,  
uniformly in $t>0,P\in M,\lambda>0$.  
The underlying idea of the proof is that  
if \eqref{introcontradfirst} holds for all $\alpha>0$,  
then for all $\alpha>0$ there exist minimizers $\ua$  
of $\Ia$, which  
approach $\{t\widetilde\xi_{P,\lambda}^g\}$ as $\alpha\to+\infty$,  
and the convergence rate is sufficiently rapid to ensure that  
for some suitable $A>0$,  
$\ua$ also satisfies \eqref{optsobolev}, uniformly in $\alpha$.  
But then $\alpha\le C$, a contradiction.\par    
In {\em Section \ref{secprelims}},  
for the reader's convenience, we establish some preliminary results   
by suitably adapting to our needs some well-known techniques  
from \cite{Tr,Au2,BN,HV,AL}.    
We show that \eqref{introcontradfirst} implies the existence  
of a minimizer $\ua\in H^1(M)$ for $\Ia$ satisfying  
$\ua\in H^1(M)$, $\ua\ge0$,  
$\intM\ua^{2^*}\dvg=1$ and such that  
\[  
\mua^{(n-2)/2}:=\max_M\ua^{-1}=:\ua(\xa)^{-1}\to0. 
\]  
We fix some small $\delta_0>0$
which
depends only on $(M,g)$.  
We show: 
\begin{align*}  
&\|\nabla_g
(\ua-\xi_{\xa,\mua^{-1}}^{g})\|_{L^2(B_{\delta_0})}
+\|\ua-\xi_{\xa,\mua^{-1}}^{g}\|_{L^{2^*}(B_{\delta_0})}\to0\\ 
&\mua^{(n-2)/2}\ua(\exp_{\xa}^{g}(\mua\ \cdot\ ))\to U 
\qquad\text{in}\ C^2_{\mathrm{loc}}(\R^n). 
\end{align*} 
The $C^2_{\mathrm{loc}}(\R^n)$-convergence 
and a change of variables imply the lower bound:  
\begin{equation}  
\label{introlowerbd}  
\|\ua\|_{L^{\br}(M,g)}\ge C^{-1}\mua^2.  
\end{equation}  
In {\em Section \ref{secuniformest}} we prove the {\em uniform estimate}:  
\[  
\ua(x)\le C\mua^{(n-2)/2}\distg(x,\xa)^{2-n}  
\qquad\forall x\in M.  
\]  
This estimate ensures a suitable decay of $\ua$   
away from $\xa$; it is a key step.
We note that pointwise estimates for minimizers to   
critical exponent equations have been   
established and used by   
Brezis and Peletier~\cite{BP},   
Atkinson and Peletier~\cite{AP},   
Rey~\cite{Rey}, Han~\cite{Ha},   
Hebey and Vaugon~\cite{HV},    
Li and Zhu~\cite{LZ,LZ3},   
Aubin and Li~\cite{AL}, and others.   
We derive our pointwise estimate    
along the line of \cite{LZ,LZ3},   
by working directly on $\ua$; new ingredients
are needed in deriving our estimate.  

In {\em Section \ref{secenergyest}},
in order to simplify calculations, we
introduce a conformal metric $\hg=\psi^{4/(n-2)}g$,
with $\psi\in C^\infty(M)$,
$\psi(\xa)=1$, $\frac 12 \le \psi\le 2$, $\|\psi\|_{ C^2}\le C$,
such that $R_{\hg}\equiv0$ in $B_{\delta_0}(\xa)$.
Our pointwise estimates in Section 2 allow us to
adapt ideas of Bahri and Coron~\cite{BC} to make an 
energy estimate of the difference:  
$\ua/\psi-\ta\xi_{\ya,\lda}^{\hg}$  
in a small ball $B_{\da}(\xa)$, 
where $\da\in[\delta_0/2,\delta_0]$, 
$\ta>0$, $\mu_\alpha^{-1}|\ya-\xa|\to 0$, $\lda>0$ are ``optimal" 
in a suitable sense.  
The main result of Section \ref{secenergyest} is the estimate  
for the projection $\ua/\psi
-\ta\xi_{\ya,\lda}^{\hg}$ on $H_0^1(B_{\da}(\xa))$,  
denoted $\wa$,  
as in Proposition \ref{propwaestimate}.\par  
In {\em Section \ref{secsplitting}} we show that by choosing  
a ``good radius" $\da\in[\delta_0/2,\delta_0]$,  
the ``boundary part" of $\ua/\psi-\ta\xi_{\ya,\lda}^{\hg}$  
may be controlled in $H^1(\partial \Ba)$,  
see Lemma \ref{lemgradbdry}.  
For $n\ge7$, the estimate resulting from
our pointwise estimates,   
Proposition \ref{propwaestimate}, Lemma \ref{lemgradbdry}  
and taking into account \eqref{epsa} is given by:  
\begin{equation}  
\label{introenergyest}  
\|\nabhg(\frac{\ua}\psi
-\ta\xi_{\ya,\lda}^{\hg})\|_{L^2(B_{\da}(\xa))}\le C\Big(\mua^2  
+(1+\mua^{-2+\beta})\alpha\|\ua\|_{L^{\br}(M,g)}^2\Big),  
\end{equation}  
where $\beta=(n-6)(n-2)/[2(n+2)]>0$ is {\em strictly positive},  
since $n\ge7$.  
By carefully exploiting orthogonality,  we  
prove the following lower bound:  
\begin{equation}  
\label{introsplitting}  
\Yg(\ua)\ge Y_{\hg}(\widetilde\xi_{\ya,\lda}^{\hg})  
+O(\mua^2\|\nabhg(\frac \ua\psi-\ta\xi_{\ya\lda}^{\hg})\|_{L^2(B_{\da}(\xa))}+\mua^{n-2}),  
\end{equation}  
see Proposition \ref{propyamabelowerbound}.\par 
At this point we have all the necessary ingredients to conclude the 
proof in the case $n\ge7$. 
We note that   
the contradiction assumption \eqref{introcontradfirst} implies:  
\begin{equation*}  
K^{-2}>\Ia(\ua)=\Yg(\ua)+\alpha\|\ua\|_{L^{\br}(M,g)}^2.  
\end{equation*}  
By the above inequality and \eqref{introsplitting}, we obtain 
\begin{align}  
\label{introalpha}  
\alpha\|\ua\|_{L^{\br}(M,g)}^2\le&K^{-2}-Y_{\hg}(\widetilde\xi_{\ya,\lda}^{\hg})\\ 
\nonumber  
&+C(\mua^2\|\nabhg(\frac \ua\psi-\ta\xi_{\ya,\lda}^{\hg})\|_{L^2(B_{\da}(\xa))}+\mua^{n-2}).  
\end{align}  
By  \eqref{aubinexpansion} (or an easy calculation since
we do not need the explicit coefficient of
$\lambda^{-4}$),  
\begin{equation}  
\label{introyamabeexp}  
|K^{-2}-Y_{\hg}(\widetilde\xi_{\ya,\lda}^{\hg})|\le C\mua^4.  
\end{equation}  
Inserting \eqref{introenergyest} and \eqref{introyamabeexp} into  
\eqref{introalpha}, and recalling that $\beta>0$, we derive  
\[  
(1+\circ(1))\alpha\|\ua\|_{L^{\br}(M,g)}^2\le C\mua^4.  
\]  
In view of \eqref{introlowerbd}, the desired contradiction $\alpha\le C$ follows, 
and Theorem \ref{thmmain} is established.\par  
Finally, {\em Section \ref{seclimitcase}}
is devoted to the proof of Theorem \ref{thmmain} in
the remaining case $n=6$.  
This is 
 more delicate than the case $n\ge 7$. Nevertheless, we can still   
obtain the   
inequality \eqref{optsobolev} with the aid of a   
uniform {\em lower} bound,  
reminiscent of an argument in \cite{LZ}.  
\paragraph{Notation}  
Henceforth, $C>0$ always denotes   
a general constant independent  
of $\alpha$, and subsequences of $\alpha\to+\infty$   
are taken without further notice.  
Denoting by $(\Omega,h)$ a Riemannian manifold  
(possibly with boundary), we set 
\begin{align*}  
&\langle\varphi,\psi\rangle_h   
=\int_\Omega\nabla_h\varphi\cdot\nabla_h\psi\,dv_h   
=\int_\Omega h^{ij}\frac{\de\varphi}{\de x^i}\frac{\de\varphi}{\de x^j}\,dv_h  
&&\forall\varphi,\psi\in H_0^1(\Omega)\\    
&\|\varphi\|_h=\sqrt{\bra\varphi,\varphi\ket_h}    
&&\forall\varphi\in H_0^1(\Omega).    
\end{align*}  
We note that the metrics $g$ and $\hg$ defined above are both equivalent to the 
Euclidean metric $\mathcal E$. 
When the specific metric is clear from the context, or irrelevant up to equivalence 
to $g$, 
we do not indicate it explicitly. 
Furthermore, for $q\ge1$ we denote: 
\begin{align*} 
&\|\varphi\|_q=\|\varphi\|_{L^q(\Ba)}\\ 
&\|U^q\|_{\stpr,\mua^{-1}}=\left(\int_{B_{\mua^{-1}}(0)}U^{\stpr q}\dy\right)^{1/\stpr}, 
\qquad\stpr=\frac{2n}{n+2}, 
\end{align*} 
where $U$ is the standard minimizer on $\R^n$ 
defined above.\par  
For ease of future reference, we prove our estimates for $n\ge3$. 
Moreover, we obtain our estimates   
for a general exponent $r\in(1,2)$,  
which could even depend on $\alpha$ (this will also be convenient for the   
local to global argument sketched in the Appendix).  
The actual value $r=\br=2n/(n+2)$ and the condition  
$n\ge6$ are used only in the final part of the proof  
of Theorem \ref{thmmain}, in Section \ref{secsplitting}  
and in Section \ref{seclimitcase}.  
\paragraph{}  
Theorem \ref{thmmain} in the case $n\ge7$ has been 
presented at the  
966th AMS Meeting at Hoboken, NJ, April 28--29, 2001.  

\section{Preliminaries}    
\label{secprelims} 
The preliminary results in this section are obtained by  
adapting standard methods to our situation,    
see, e.g., \cite{Tr,Au2,BN,HV,AL}.  
For the reader's convenience, we sketch their proofs.
Throughout this section, we assume $n\ge3$.\par    
For every $\alpha>0$ and for $r\in(1,2)$ (possibly depending on $\alpha$) 
we consider the functional: 
\begin{equation*}    
\Ia(u)    
=\frac{\intM\{|\nabg u|^2+c(n)\Rg u^2\}\dvg    
+\alpha\|u\|_{L^r(M)}^2}    
{\|u\|_{L^{2^*}(M)}^2},    
\end{equation*}    
defined for all $u\in H^1(M)\setminus\{0\}$.    
If \eqref{optsobolev} is false, then for all $\alpha>0$    
we have    
\begin{equation}    
\label{contradassumption}    
\inf_{H^1(M)\setminus\{0\}}I_\alpha<K^{-2}.    
\end{equation}    
\begin{Proposition}[Existence of a minimizer]    
\label{propminimizerexist}    
For all $\alpha>0$ there exists a non-negative
minimizer $\ua\in H^1(M)$ such that    
\begin{align*}    
&I_\alpha(\ua)=\la=\inf_{H^1(M)\setminus\{0\}}I_\alpha<K^{-2}\\    
&\intM\ua^{2^*}\dvg=1.    
\end{align*}    
Moreover, $\ua\in C^{2,r-1}(M)$ is a classical solution of the 
Euler-Lagrange equation:    
\begin{equation}    
\label{uaeq}    
-\Delta_g\ua    
+c(n)\Rg\ua+\alpha\|\ua\|_{L^r(M)}^{2-r}\ua^{r-1}    
=\la\ua^{2^*-1}\qquad\text{on}\ M.   
\end{equation}    
\end{Proposition}    
\begin{proof}    
By homogeneity, it is equivalent to minimize $I_\alpha$ on the set    
\[    
\mathcal A=\{u\in H^1(M):\ \intM|u|^{2^*}\dvg=1\}.    
\]    
However, $\mathcal A$ is not sequentially weakly closed in $H^1(M)$. Therefore,    
as usual,    
 for fixed $\alpha$ and     
for all $1\le q<2^*$ we define:    
\[    
\Aq=\{u\in H^1(M):\ \intM|u|^q\dvg=1\}    
\]    
and we consider the functional    
\[    
I_q(u)=\frac{\intM\{|\nabg u|^2+c(n)\Rg u^2\}\dvg    
+\alpha\big(\intM|u|^r\dvg\big)^{2/r}}    
{\big(\intM|u|^q\dvg\big)^{2/q}}   
\]    
on $\Aq$.     
By standard arguments $\inf_{\Aq}I_q$ is attained, i.e.,     
for every $1\le q<2^*$ there exists    
$\uq\in\Aq$ such that    
\[    
I_q(\uq)=\inf_{\Aq}I_q=:\ell_q.    
\]    
The minimizer $\uq$ satisfies the Euler-Lagrange equation:    
\begin{equation}    
\label{uqeq}    
-\Delg\uq+c(n)\Rg\uq+\alpha\|\uq\|_{L^r(M)}^{2-r}\uq^{r-1}=\ell_q\uq^{q-1} 
\qquad\text{on}\ M.    
\end{equation}    
The sequence $\uq$ is bounded in $H^1(M)$; therefore passing to a subsequence    
we can assume that there exists $\ua\in H^1(M)$ such that    
$\uq\weak\ua$ weakly in $H^1(M)$,    
strongly in $L^2(M)$ and a.e.    
Since for every fixed $u$ we have $I_q(u)\to I_\alpha(u)$ as $q\to2^*$,    
it is clear that     
\[    
\limsup_{q\to\infty}\ell_q\le\la<K^{-2}.    
\]    
Consequently, for every $0<2^*-q\ll1$,    
we can apply the Moser iteration technique 
to \eqref{uqeq} to derive
a uniform bound $\sup_M\uq\le C(\alpha)$, 
where $C(\alpha)>0$ is a constant independent of $q$  
(see, e.g.,  \cite{AL}).    
Then by dominated convergence, $\ua\in\mathcal A$    
and by weak semicontinuity $I_\alpha(\ua)\le\liminf\ell_q\le\la$.    
The $\ua$ is a desired minimizer.\par
The proof of the existence of the minimizer shows that $\ua$    
is in $L^\infty(M)$ for every {\em fixed} $\alpha$.    
Then standard elliptic theory implies that $\ua\in C^{1,\beta}(M)$    
for some $0<\beta<1$. Therefore $\ua^{r-1}\in C^{0,r-1}(M)$,     
and by Schauder estimates $\ua\in C^{2,r-1}(M)$.    
\end{proof}    
\begin{Remark}    
Since $0<r-1<1$, the nonlinearity $u^{r-1}$ is {\em sublinear}    
and therefore we can not use the maximum principle to conclude $\ua>0$ on $M$.    
\end{Remark}    
\begin{Proposition}[Standard blowup]    
\label{propstandblowupua}    
As $\alpha\to+\infty$, we have:    
\begin{align*}    
\tag{i}    
&\ua\to0,\quad\text{weakly in}\;H^1(M), 
\ \text{strongly in}\ L^p(M)\ \forall1\le p<2^*    
\ \text{and}\;a.e.\\    
\tag{ii}    
&\intM|\nabg\ua|^2\dvg\to K^{-2}\\    
\tag{iii}    
&\alpha\|\ua\|_{L^r(M)}^2\to0\\    
\tag{iv}    
&\la\to K^{-2}\\   
\tag{v}   
&\max_M\ua\to+\infty.    
\end{align*}    
\end{Proposition}    
\begin{proof}    
By compactness,    
for any $\eps>0$ there exists $C_\eps>0$ such that:    
\[    
\|u\|_{L^2(M)}^2\le\eps\intM|\nabg u|^2\dvg 
+C_\eps\|u\|_{L^r(M)}^2.    
\]    
So,    
\begin{align*}    
I_\alpha(\ua)=&\la=\intM\big\{|\nabg\ua|^2+c(n)\Rg\ua^2\big\}\dvg    
+\alpha\|\ua\|_{L^r(M)}^2\\    
\ge&(1-\eps\,c(n)\max_M|\Rg|)\intM|\nabg\ua|^2\dvg    
+(\alpha-C_\eps)\|\ua\|_{L^r(M)}^2.    
\end{align*}    
Fixing a small $\eps$ we obtain:    
\[    
\frac{1}{2}\intM|\nabg\ua|^2\dvg+(\alpha-C_\eps)\|\ua\|_{L^r(M)}^2    
\le I_\alpha(\ua)<K^{-2}.    
\]    
Consequently,  
\[    
\intM|\nabg\ua|^2\dvg\le C,\qquad\alpha\|\ua\|_{L^r(M)}^2\le C,    
\]    
and therefore, 
\[    
\intM\ua^r\dvg\to0\qquad\text{as}\;\alpha\to+\infty.    
\]    
Passing to a subsequence, we have (i).    
Furthermore, we can assume that for some $\theta,\eta\in[0,+\infty)$ there 
holds (along a subsequence):
\[    
\intM|\nabg\ua|^2\dvg\to\theta\quad\text{and}\quad    
\alpha\|\ua\|_{L^r(M)}^2\to\eta,\quad \text{as}\;\alpha\to+\infty.
\]    
Proof of (ii)--(iii). We have to show that $\theta=K^{-2}$ and $\eta=0$.    
By the Sobolev inequality as in \cite{Au}, for every $\eps>0$
there exists $A_\eps>0$ such that:    
\begin{equation}    
\label{aubinsobolev}    
\|\ua\|_{L^{2^*}(M)}^2    
\le K^2(1+\eps)\intM|\nabg\ua|^2\dvg+A_\eps\|\ua\|_{L^r(M)}^2.    
\end{equation}    
Letting $\alpha\to+\infty$ in \eqref{aubinsobolev} we obtain:    
\[    
1=\big(\intM\ua^{2^*}\dvg\big)^{2/2^*}\le K^2(1+\eps)\,\theta.    
\]    
Sending $\eps\to0$, we conclude $1\le K^2\theta$.    
On the other hand, we have by definition of $\ua$:    
\[    
\intM\{|\nabg\ua|^2+c(n)\Rg\ua^2\}\dvg    
+\alpha\|\ua\|_{L^r(M)}^2=\la<K^{-2}.    
\]    
Sending $\alpha\to+\infty$ we find $\theta+\eta\le K^{-2}$.    
It follows that $\theta=K^{-2}$ and $\eta=0$, as asserted.\\   
Proof of (iv). This is an immediate consequence of (i)--(ii)--(iii)   
and the definition of $\la$.\\   
Proof of (v).    
We have:    
\begin{align*}    
1=\intM\ua^{2^*}\dvg\le(\max_M\ua)^{2^*-r}\intM\ua^r\dvg    
=\circ(1)(\max_M\ua)^{2^*-r}.    
\end{align*}    
\end{proof}    
Our next aim is to show that, after rescaling, the limit profile of $\ua$    
is the standard minimizer $U$, and that $\ua$ approaches this limit     
``in energy", as in Proposition \ref{propstandblowupva}.\\    
Let $\xa\in M$ be a maximum point of $\ua$, namely    
$\ua(\xa)=\max_M\ua$, then by Proposition \ref{propstandblowupua}--(v)   
we have    
\begin{equation}    
\label{muaeq}    
\mua:=\ua(\xa)^{-2/(n-2)}\to0\quad\text{as}\;\alpha\to+\infty.    
\end{equation}    
Let $\delta_0>0$ be a small constant to be fixed below
(e.g., less than injectivity radius).    
Let $\delta_0/2\le\delta_\alpha\le\delta_0$.

\begin{Proposition}[Convergence in energy]\ \
\label{propstandblowupva}
\begin{equation}
\label{vastrongconv}
\lim_{\alpha\to+\infty}\int_{ B_{\delta_\alpha}(\xa) }
\big\{|\nabla_g(\ua-\xi^g_{\xa,\mu_\alpha^{-1}})|^2
+|\ua-\xi^g_{\xa,\mu_\alpha^{-1}}|^{2^*}\big\}dv_g=0.
\end{equation}
\end{Proposition}

\begin{proof}
We consider the following rescaling of $\ua$ on
the geodesic ball $B_{\delta_\alpha}(\xa)$:
\begin{equation}
\va(y)
=\mua^{(n-2)/2}\ua(\exp_{\xa}(\mua y)),\qquad y\in\Oa,
\label{100}
\end{equation}
where
\begin{equation}
\Oa=\mua^{-1}\exp_{\xa}^{-1}(B_{\delta_\alpha}(\xa))
=\mua^{-1}B_{\delta_\alpha}(0).
\label{200}
\end{equation}
$\va$ satisfies    
\begin{align}    
\label{vaeq}    
&-\Delga\va+c(n)R_{g_\alpha}v_\alpha+
\epsa\va^{r-1}=\la\va^{2^*-1}    
\quad\text{in}\;\Oa,   
\end{align}    
where 
$$ 
\ga(y)=g(\exp_{\xa}(\mua y)),\quad
|R_{g_\alpha}|\le C\mu_\alpha^2,
$$
and 
\[    
\epsa:=\alpha\mua^{n-\frac{n-2}{2}r}\|\ua\|_{L^r(M)}^{2-r}.    
\]    
We observe that the rescaled  metric $\ga$ converges     
to the Euclidean metric $(\delta_{ij})$ on $\R^n$ uniformly on compact subsets,    
and it is equivalent to $(\delta_{ij})$, uniformly in $\alpha$,    
i.e., there exists $C>0$ independent of $\alpha$ such that    
$C^{-1}\delta_{ij}\le{\ga}_{,ij}(y)\le C\delta_{ij}$.    
We claim that    
\begin{equation}    
\label{epsa}    
\epsa\le\alpha\|\ua\|_{L^r(M)}^2\to0    
\qquad\text{as}\;\alpha\to+\infty.    
\end{equation}    
Indeed, by the definition of $\mua$ and $\epsa$,    
\begin{align*}    
\epsa=\frac{\alpha\|\ua\|_{L^r(M)}^2}{(\max_M\ua)^{2^*-r}\intM\ua^r\dvg}    
\end{align*}    
and    
\[    
1=\intM\ua^{2^*}\dvg\le(\max_M\ua)^{2^*-r}\intM\ua^r\dvg.    
\]    
Property \eqref{epsa} now follows by Proposition \ref{propstandblowupua}--(iii).   
By a change of variables,     
\[    
\inta\va^{2^*}\dvga=\int_{\Ba}\ua^{2^*}\dvg    
\]    
and  
\begin{align*}   
\inta|\nabga\va|^2\dvga=& 
\int_{\Ba}\{|\nabg\ua|^2+c(n)\Rg\ua^2\}\dvg. 
\end{align*}   
Consequently, by the definition of $\ua$
\begin{equation}    
\label{limsupva2star}    
\limsup_{\alpha\to+\infty}\inta\va^{2^*}\dvga\le1    
\end{equation}    
and  by Proposition \ref{propstandblowupua}--(i)--(ii),
\begin{equation}    
\label{limsupgradva}    
\limsup_{\alpha\to+\infty}\inta|\nabga\va|^2\dvga\le K^{-2}.    
\end{equation}    
By the  definition of $\mua$,
$\va(y)\le\va(0)=1$, thus, by
 standard elliptic estimates, there exists $v\in C_{\text{loc}}^1(\R^n)$    
such that, along a subsequence, $\va\to v$ in $C_{\text{loc}}^1(\R^n)$,
and $v(0)=1$. 
Furthermore, $v$ satisfies:    
\begin{equation}    
\label{nablav}    
\intRn|\nabla v|^2\dy    
=\lim_{R\to+\infty}\int_{B_R}|\nabla v|^2\dy    
=\lim_{R\to+\infty}\lim_{\alpha\to+\infty}\int_{B_R}|\nabga\va|^2\dvga    
\le K^{-2},    
\end{equation}    
and    
\begin{equation}    
\label{v2star}    
\intRn v^{2^*}\dy=\lim_{R\to+\infty}\int_{B_R}v^{2^*}\dy    
=\lim_{R\to+\infty}\lim_{\alpha\to+\infty}\int_{B_R}\va^{2^*}\dvga\le1.    
\end{equation}    
In particular, $v\in D^{1,2}(\R^n)$, and taking pointwise limits in \eqref{vaeq}    
we find that $v$ satisfies:    
\begin{align*}    
&-\Delta v=K^{-2}v^{2^*-1}\qquad\qquad\text{in}\;\R^n\\    
&0\le v\le 1,\ v(0)=1.    
\end{align*}    
Multiplying the above equation by $v$ and integrating by parts,    
and recalling the definition of $K$ we have:    
\[    
K^{-2}\intRn v^{2^*}\dy=\intRn|\nabla v|^2\dy    
\ge K^{-2}\big(\intRn v^{2^*}\dy\big)^{2/2^*}.    
\]    
Therefore, $\big(\intRn v^{2^*}\dy\big)^{1-2/2^*}\ge1$,    
which together with \eqref{nablav} and \eqref{v2star}    
implies $\intRn v^{2^*}\dy=1$ and $\intRn|\nabla v|^2\dy=K^{-2}$    
and thus necessarily $v=U$.   
Since the limit $v$ is independent of
subsequences, the convergence is for all
$\alpha\to+\infty$ with $x_\alpha\to P$.
At this point, it is intuitively clear that  
Proposition \ref{propstandblowupua}--(ii) should imply the 
``strong convergence" \eqref{vastrongconv}; however  
we face some minor technicality 
due to the fact that $\va$ does not 
necessarily vanish on $\de\Oa$. 
Using the elementary calculus inequality:    
\[    
\big||a+b|^p-|a|^p-|b|^p\big|\le C(p)\big(|a|^{p-1}|b|+|a||b|^{p-1}\big),    
\qquad\forall a,b\in\R^n, p\ge1    
\]    
with $p=2^*$, $a=U$, and $b=\va-U$,    
we have:    
\begin{align*}    
\inta|\va-U|^{2^*}&\dvga\le\inta\va^{2^*}\dvga-\inta U^{2^*}\dvga\\    
&+C\big(\inta U^{2^*-1}|\va-U|\dvga+\inta U|\va-U|^{2^*-1}\dvga\big)\\ 
\le&\circ(1)+C\big(\inta
U^{2^*-1}|\va-U|\dvga+\inta U|\va-U|^{2^*-1}\dvga\big).    
\end{align*}    
The right hand side is easily seen to vanish    
as $\alpha\to+\infty$:    
\begin{align*}    
\inta U|\va-&U|^{2^*-1}\dvga\\   
=&\int_{B_R}U
|\va-U|^{2^*-1}\dvga+\int_{\Oa\setminus B_R} U|\va-U|^{2^*-1}\dvga\\    
\le&\int_{B_R} U|\va-U|^{2^*-1}\dvga\\    
&\qquad\qquad+\big(\int_{\Oa\setminus B_R}U^{2^*}\dvga\big)^{1/2^*}    
\big(\int_{\Oa\setminus B_R}|\va-U|^{2^*}\dvga\big)^{1/{2^*}'}\\    
\le&\int_{B_R}U|\va-U|^{2^*-1}\dvga    
+C\int_{\R^n\setminus B_R}U^{2^*}\dy.    
\end{align*}    
By taking $R$ large, the second integral can be made arbitrarily small;    
then, by $C_{\text{loc}}^1$-convergence, the first integral    
is small for large $\alpha$.    
Hence,    
\[    
\lim_{\alpha\to+\infty}\inta U|\va-U|^{2^*-1}\dvga=0.    
\]    
Similarly, one easily checks that    
\[    
\lim_{\alpha\to+\infty}\inta U^{2^*-1}|\va-U|\dvga=0.    
\]    
The strong convergence of the gradients is straightforward:    
\begin{align*}    
\big|\inta&\nabga(\va-U)\cdot\nabga U\dvga\big|    
\le\int_{B_R}|\nabga(\va-U)||\nabga U|\dvga\\    
&+\big(\int_{\Oa\setminus B_R}|\nabga(\va-U)|^2\dvga\big)^{1/2}    
\big(\int_{\Oa\setminus B_R}|\nabga U|^2\dvga\big)^{1/2}\\    
\le&\int_{B_R}|\nabga(\va- U)||\nabga U|\dvga    
+C\big(\int_{\Oa\setminus B_R}|\nabla U|^2\dy\big)^{1/2}    
\end{align*}    
and therefore    
\[    
\lim_{\alpha\to+\infty}\inta\nabga(\va-U)\cdot\nabga U\dvga=0.    
\]    
Consequently, by \eqref{limsupgradva} and  
since $\inta  |\nabga U|^2\dvga\to K^{-2}$, we conclude:    
\begin{align*}    
\inta|\nabga(\va-U)&|^2\dvga    
=\inta|\nabga\va|^2\dvga-\inta|\nabga U|^2\dvga\\    
&-2\inta\nabga(\va-U)\cdot\nabga U\dvga    
\le o(1),    
\end{align*}    
and \eqref{vastrongconv} follows after a change of variables.    
\end{proof}    
\begin{Corollary}[One point concentration for $\ua$]    
\label{coronepointconc}  
For  any $\eps>0$ there exist $\delta_\eps>0$
and $\alpha_\eps>0$ such that
\begin{equation*}
\int_{M\setminus B_{\mua/\delta_\eps}(\xa)}
\big\{|\nabg\ua|^2+\ua^{2^*}\big\}\dvg\le\eps
\end{equation*}
for all $\alpha\ge\alpha_\eps$.
In particular, for any fixed $\rho>0$,   
\begin{equation*}    
\lim_{\alpha\to+\infty}\int_{M\setminus B_\rho(\xa)} 
\big\{|\nabg\ua|^2+\ua^{2^*}\big\}\dvg=0.    
\end{equation*}    
\end{Corollary}    
\begin{proof} For any $\eps>0$, by (\ref{vastrongconv}) and
a change of variable, there exists $\delta_\eps>0$ and $\alpha_\eps'$
such that for all $\alpha\ge \alpha_\eps'$,
$$
\int_{ B_{\mu_\alpha/\delta_\eps}(x_\alpha) }
|\nabla_g u_\alpha|^2dv_g\ge
\int_{\R^n}|\nabla U|^2-\frac \eps 4=K^{-2}-\frac \eps 4,
$$
and
$$
\int_{ B_{\mu_\alpha/\delta_\eps}(x_\alpha) }
u_\alpha^{ 2^*}dv_g\ge
\int_{\R^n} U^{ 2^*}-\frac \eps 4=1-\frac \eps 4.
$$
Recall that $\int_M|\nabla_g u_\alpha|^2dv_g\to K^{-2}$ and
$\int_Mu_\alpha^{ 2^*}dv_g=1$, we can take some 
$\alpha_\eps\ge \alpha_\eps'$ 
such that for all $\alpha\ge \alpha_\eps$, 
$$
\int_{M\setminus B_{\mua/\delta_\eps}(\xa)}
\big\{|\nabg\ua|^2+\ua^{2^*}\big\}\dvg\le\eps.
$$    
\end{proof} 
\begin{Corollary} 
\label{coreasysupest} 
For any fixed $\rho>0$,  
\[ 
\lim_{\alpha\to+\infty}\|\ua\|_{L^\infty(M\setminus B_\rho(\xa))}=0. 
\] 
\end{Corollary}  
\begin{proof} 
Equation \eqref{uaeq} implies the differential inequality 
\[ 
-\Delg\ua+[c(n)\Rg-\la\ua^{2^*-2}]\,\ua\le0 
\qquad\text{on}\ M. 
\] 
By Corollary \ref{coronepointconc} we have, 
for all $x\in M\setminus B_{\rho}(\xa)$, that
\newline 
$\|\ua^{2^*-2}\|_{L^{n/2}(B_{\rho/2}(x))} 
=\|\ua\|_{L^{2^*}(B_{\rho/2}(x))}^{2^*-2}
\le \|\ua\|_{L^{2^*}(M\setminus B_{\rho/2}(x_\alpha)) }=\circ(1)$. 
By Moser iterations, we derive 
\[ 
\|\ua\|_{L^\infty(B_{\rho/4}(x))}\le
C\|\ua\|_{L^1(B_{\rho/2}(x))}\le  C\|\ua\|_{L^1(M)}, 
\] 
and the claim follows by Proposition \ref{propstandblowupua}--(i). 
\end{proof} 


\section{Uniform estimate}    
\label{secuniformest} 
The $C_{\text{loc}}^2(\R^n)$-convergence of     
the rescaled minimizer $\va$ to $U$ readily    
provides a complete description of $\ua$ in a ball of  
{\em shrinking} radius $B_{\rho\mua}(\xa)$,    
for any $\rho>0$.    
In particular, it implies the estimate:    
\[
\ua(x)=\mua^{-(n-2)/2}\va\big(\mua^{-1}(\exp_{\xa}g)^{-1}(x)\big)
\le C(\rho)\mua^{-(n-2)/2},    
\ \ \forall x\in B_{\rho\mua}(\xa),    
\]
and consequently    
\begin{equation}    
\label{estimatenearxa}    
\ua(x)\le C(\rho)\mua^{(n-2)/2}\distg(x,\xa)^{2-n}    
\qquad\forall x\in B_{\rho\mua}(\xa).    
\end{equation}    
Our aim in this section is to show that \eqref{estimatenearxa}    
holds {\em uniformly on} $M$. This type of estimate for minimizers
has been obtained by
Brezis and Peletier~\cite{BP}
and by Atkinson and Peletier~\cite{AP} in the
radially symmetric case 
on Euclidean balls, by
Rey~\cite{Rey} and
 Han~\cite{Ha} on general domains
in $\R^n$,    
and by Hebey and Vaugon~\cite{HV},
Li and Zhu~\cite{LZ} and  Aubin and Li~\cite{AL}
on Riemannian manifolds.    
Our approach, similar in spirit to \cite{LZ}, requires new ingredients. 
Throughout this section, we assume $n\ge3$. 
\begin{Proposition}    
\label{propuasupest}    
For every $\alpha$ sufficiently large, $\ua$ satisfies    
\begin{equation}
\ua(x)\le C\mua^{(n-2)/2}\distg(x,\xa)^{2-n}\quad\forall x\in M.    
\label{unifest}
\end{equation}
Here $C>0$ is a constant depending on $(M,g)$ only.    
Consequently, we have the following {\em uniform estimate}   
for $\va$:   
\begin{align}   
\label{uniformest}   
\va(y)\le\frac{C}{1+|y|^{n-2}}, \qquad\forall y\in\Oa.   
\end{align}   
\end{Proposition}    
We shall prove Proposition \ref{propuasupest} by showing that    
\begin{equation}    
\label{ualephia}    
\ua(x)\le C\phia(x)\quad\forall x\in M,    
\end{equation}    
for some $\phia>0$     
satisfying:    
\begin{equation}    
\label{phiablowuprate}    
C^{-1}\mua^{(n-2)/2}\distg(x,\xa)^{2-n}\le\phia(x)\le C\mua^{(n-2)/2}\distg(x,\xa)^{2-n}    
\quad\forall x\in M    
\end{equation}    
for every $\alpha$ sufficiently large.    
In fact, our main effort will be  to construct     
a suitable such  $\phia$.    
We set    
\[    
\zea=\frac{\ua}{\phia}.    
\]    
We have to show $\zea\le C$ pointwise on $M$.    
By the conformal invariance,  $\zea$  satisfies   
\begin{align}    
\label{zeaeq}    
&-\Delgh\zea\\  
\nonumber
=&\la\zea^{2^*-1}    
-\phia^{1-2^*}(-\zea\Delta_g\phia+c(n)\Rg\ua+\alpha\|\ua\|_{L^r(M)}^{2-r}\ua^{r-1})    
\;\text{in}\;M\setminus\{\xa\},    
\end{align}    
where $\gh$ is the metric conformal to $g$ defined     
in terms of $\phia$ by $\gh=\phia^{4/(n-2)}g$.\par    
Indeed, we have 
\[   
-\Delta_{\gh}\frac{u}{\phia}+c(n)R_{\gh}\frac{u}{\phia} 
=\phia^{1-2^*}\big(-\Delg u+c(n)\Rg u\big),   
\quad \forall\ u\in C^2(M\setminus\{x_\alpha\}).   
\]   
Taking $u=\phia$, we obtain   
\begin{equation*}   
c(n)R_{\gh}=\phia^{1-2^*}(-\Delg\phia+c(n)\Rg\phia).   
\end{equation*}   
Taking    
$u=\ua$,  we find    
\begin{equation*}    
-\Delgh\zea=-\phia^{1-2^*}\Delg\ua    
-c(n)(R_{\gh}-\phia^{-4/(n-2)}\Rg)\zea.    
\end{equation*}    
It follows that    
\begin{align*}    
-\Delgh\zea=\phia^{-(n+2)/(n-2)}(-c(n)\Rg\ua 
-\alpha\|\ua\|_{L^r(M)}^{2-r}\ua^{r-1}+&\la\ua^{2^*-1)})\\    
+&\phia^{1-2^*}\zea\Delg\phia,    
\end{align*}    
which implies \eqref{zeaeq}.\par    
By the uniform estimate \eqref{phiablowuprate},    
the metrics $\gh$ satisfy a Sobolev inequality with a constant    
{\em independent} of $\alpha$:    
\begin{Lemma} 
\label{lemghatsobolev}    
There exists a constant $C>0$ independent of $\alpha$    
such that for all $u\in H^1(M)$, $u\equiv0$    
in a neighborhood of $x_\alpha$:    
\begin{equation}    
\left(\intM|u|^{2^*}\dvgh\right)^{2/2^*}    
\le C\intM|\nabgh u|^2\dvgh.    
\end{equation}    
\end{Lemma}    
\begin{proof}    
It is well-known 
(see, e.g.,  Appendix A in \cite{LZ})    
that there exists a constant $C=C(M,g)$ such that    
for all $x_0\in M$, $u\in H^1(M)$, $u\equiv0$ in     
a neighborhood of  $x_0$, there holds:    
\begin{equation}    
\label{weightedsob}    
\left(\intM\frac{|u|^{2^*}}{\distg(x,x_0)^{2n}}\dvg\right)^{2/2^*}    
\le C\intM\frac{|\nabg u|^2}{\distg(x,x_0)^{2n-4}}\dvg.    
\end{equation}    
Now it suffices to observe that by conformality of $\gh$ we have:    
\[    
\dvgh=\phia^{2^*}\dvg\qquad\text{and}\qquad    
|\nabgh u|^2=\phia^{-4/(n-2)}|\nabg u|^2,    
\]    
and to recall \eqref{phiablowuprate}.    
\end{proof}    
At this point it is clear from \eqref{zeaeq} that if we can    
find a function    
$\phia>0$ satisfying \eqref{phiablowuprate}    
and such that:    
\begin{equation}    
\label{zeaineq}    
-\frac{\ua}{\phia}\Delg\phia+c(n)\Rg\ua+\alpha\|\ua\|_{L^r(M)}^{2-r}\ua^{r-1}\ge0    
\quad\text{in}\;M\setminus\{\xa\},    
\end{equation}    
then the corresponding $\zea$ will satisfy:
\begin{align}    
\label{zeapb}    
&-\Delgh\zea\le\la\zea^{2^*-1}\quad\text{in}\;M\setminus\{\xa\}\\    
\nonumber 
&\int_{M\setminus B_{\mua/\delta_1}(\xa)}\zea^{2^*}\dvgh    
=\int_{M\setminus B_{\mua/\delta_1}(\xa)}\ua^{2^*}\dvg\le\eps,    
\end{align}    
(recall Corollary \ref{coronepointconc} in Section \ref{secprelims}).    
For any $\rho>0$, let    
\begin{equation*}    
R_i:= \frac{[2-2^{-(i-1)}]\mua}{\rho},\qquad i=1,2,3,\dots    
\end{equation*}                                                   
By \eqref{phiablowuprate} we may choose    
cutoff functions $\eta_i$ (depending on $\alpha$)    
satisfying:    
\begin{align*}    
&\eta_i\equiv1\qquad\text{in}\;M\setminus B_{R_{i+1}}\\    
&\eta_i\equiv0\qquad\text{in}\;M\setminus B_{R_i}\\    
&|\nabgh\eta_i|\le C(\rho)2^i,\quad|\nabgh^2\eta_i|\le C(\rho)4^i.    
\end{align*}                                                      
Then we shall have all necessary ingredients to apply    
the Moser iteration technique to \eqref{zeapb} and to derive:   
\begin{Lemma}    
\label{lemzeaest}    
The following pointwise upper bound holds:   
\begin{equation}    
\label{zeaestclaim}    
\zea\le C\qquad\text{in}\;M\setminus B_{\mua/\delta_0}(\xa).    
\end{equation}    
\end{Lemma}    
\begin{proof} 
By applying Moser iterations to \eqref{zeapb}, see \cite{LZ} for the detailed proof. 
\end{proof}    
Estimates \eqref{estimatenearxa} and \eqref{zeaestclaim}   
will then imply \eqref{ualephia} and thus     
Proposition \ref{propuasupest} will be established.\par   
We note that \eqref{zeaineq} is trivially satisfied if $\ua=0$.    
In $(M\setminus\{\xa\})\cap\{\ua>0\}$, \eqref{zeaineq} is equivalent to:    
\begin{equation}    
\label{phiaineq}    
-\Delg\phia+    
\left[c(n)\Rg+\alpha\left(\frac{\|\ua\|_{L^r(M)}}{\ua}\right)^{2-r}\right]\phia\ge0,    
\end{equation}    
and the operator on the left hand side above is linear in $\phia$.    
Furthermore, the blowup rate as in \eqref{phiablowuprate} is satisfied    
if $\mua^{(2-n)/2}\phia$ has the blowup rate of the Green's function with pole at $\xa$.    
In fact, we shall obtain a $\phia$ of the form $\phia=\mua^{(n-2)/2}\Ga$,    
with $\Ga$ the Green's function for the operator $-\Delg+\Vatilde$    
with pole at $\xa$,    
and where $\Vatilde$ is a truncation of the ``potential"    
$c(n)\Rg+\alpha(\|\ua\|_{L^r(M)}/\ua)^{2-r}$ appearing in \eqref{phiaineq}.    
The detailed proof follows.\par   
We define a function $\Vatilde$ in the following way:    
\begin{equation*}    
\Vatilde:=    
\begin{cases}    
\min\left\{c(n)\Rg+\alpha\left(\frac{\|\ua\|_{L^r(M)}}{\ua}\right)^{2-r},1\right\}    
\quad&\text{if}\;\ua\neq0\\    
1&\text{if}\;\ua=0.    
\end{cases}    
\end{equation*}    
Note that $\Vatilde$ is Lipschitz on $M$     
(with Lipschitz constant depending on $\alpha$) and it is uniformly bounded:    
\begin{equation}    
\label{Vatildelinftyest}    
-c(n)\|\Rg\|_\infty\le\Vatilde\le1.    
\end{equation}    
We shall prove \eqref{ualephia} with $\phia=\mua^{(n-2)/2}\Ga$    
and $\Ga$ defined in the following    
\begin{Proposition}    
\label{propGa}    
The operators $-\Delg+\Vatilde$ are coercive on $H^1(M)$ for sufficiently    
large  $\alpha$, with coercivity constant uniform in      
$\alpha$. Consequently,    
for every $\alpha$ sufficiently large there exists a unique (distributional)    
solution $\Ga$ to the equation:    
\begin{equation}    
\label{Gaeq}    
-\Delta_g\Ga+\Vatilde\Ga=\delxa,\qquad\text{on}\;M.    
\end{equation}    
Furthermore, the first nonzero eigenvalue of $-\Delg+\Vatilde$    
is bounded away from zero and therefore    
$\Ga$ satisfies, for some constant $C>0$   
independent of $\alpha$,  
\begin{itemize}    
\item[(i)]    
$\Ga\in C_{\text{loc}}^2(M\setminus\{\xa\})$;  
\item[(ii)]   
$C^{-1}\distg(x,\xa)^{2-n}\le\Ga(x)\le C\distg(x,\xa)^{2-n}  
\quad\forall\ x\in M$;  
\item[(iii)]   
$\ua\Delta_g\Ga\le[c(n)\Rg\ua+\alpha\|\ua\|_r^{2-r}\ua^{r-1}]\Ga$   
in $M\setminus\{\xa\}$.  
\end{itemize}    
\end{Proposition}    
In order to prove Proposition \ref{propGa} we need the following    
\begin{Lemma}    
\label{lemVcoercivity}   
The functions $\Vatilde$ satisfy:        
\[   
\lim_{\alpha\to+\infty}\volg\{\Vatilde<\frac 12\}=0.    
\]    
\end{Lemma}     
\begin{proof}      
Note that for every measurable set $E$ such that    
$\overline E\subset M\cap\{\ua>0\}$     
we have the lower bound:    
\[   
\|\ua\|_{L^r(E)}\|\ua^{-1}\|_{L^r(E)}\ge(\volg E)^{2/r}.    
\]   
Indeed, using the H\"older inequality we find:    
\begin{align*}    
\volg E=&\int_E\,dv_g=\int_E\ua^{r/2}\ua^{-r/2}\dvg    
\le\|\ua\|_{L^r(E)}^{r/2}\|\ua^{-1}\|_{L^r(E)}^{r/2} .   
\end{align*}    
It follows that   
\begin{align}   
\label{lowerest}   
\|(\|\ua\|_{L^r(M)}\ua^{-1})^{2-r}\|_{L^{r/(2-r)}(E)}   
=&\|\ua\|_{L^r(M)}^{2-r}\|\ua^{-(2-r)}\|_{L^{r/(2-r)}(E)}\\   
\nonumber   
\ge&\|\ua\|_{L^r(E)}^{2-r}\|\ua^{-1}\|_{L^r(E)}^{2-r}\ge|E|^{(2-r)2/r}.   
\end{align}   
Let $E_\alpha:=\{\widetilde V_\alpha<1/2\}$.     
Then $\overline E_\alpha\subset M\cap\{\ua>0\}$ and therefore,   
by (\ref{lowerest}),   
\[   
(\volg E_\alpha)^{(2-r)2/r}\le   
\|(\|\ua\|_{L^r(M)}\ua^{-1})^{2-r}\|_{L^{r/(2-r)}(E_\alpha)}.   
\]   
On the other hand, since    
\[   
\alpha(\|u_{\alpha}\|_{L^r(M)}u_{\alpha}^{-1})^{2-r}   
<\frac 12 +c(n)|R_g|,\qquad   
\mbox{on}\ E_\alpha,   
\]   
we have   
\[   
\alpha\|(\|\ua\|_{L^r(M)}\ua^{-1})^{2-r}\|_{L^{r/(2-r)}(E_\alpha)}   
\le (\frac 12 +c(n)\|R_g\|_{L^\infty(M)})(\volg M)^{(2-r)/r},   
\]   
and consequently,   
\[   
\alpha(\volg E_\alpha)^{(2-r)/r}\le C,   
\]   
for some $C>0$ independent of $\alpha$.   
Now Lemma \ref{lemVcoercivity} follows immediately.   
\end{proof}    
\begin{proof}[Proof of Proposition \ref{propGa}]    
Proof of the coercivity.    
For $\tilde\gamma=1/2$ and      
$u\in H^1(M)$, by the Sobolev inequality and     
a straightforward computation we have:    
\begin{align*}    
\intM\{&|\nabg u|^2+    
\Vatilde u^2\}\dvg=\intM\{|\nabg u|^2+\tilde\gamma    
u^2+(\Vatilde -\tilde\gamma )u^2\}\dvg\\    
\ge&\intM\{|\nabg u|^2+\tilde\gamma u^2-(\Vatilde     
-\tilde\gamma)_-u^2\}\dvg\\    
\ge&\intM\{|\nabg u|^2+\tilde\gamma    
u^2\}\dvg-\|(\Vatilde-\tilde\gamma)_-\|_{L^{n/2}(M)}\|u\|_{L^{2^*}(M)}^2\\    
\ge&\intM\{|\nabg u|^2+\tilde\gamma u^2\}\dvg    
-C\volg\{\Vatilde<1/2\} \intM\{|\nabg u|^2+u^2\}\dvg,   
\end{align*}    
where $(\Vatilde-\tilde\gamma)_- \ge 0$ denotes the negative part 
of $\Vatilde-\tilde\gamma$.   
The coercivity and its uniformity in $\alpha$    
follow from the above and Lemma \ref{lemVcoercivity}.\\   
Proof of (i) and (ii).  Because of the   
coercivity of $-\Delg+\Vatilde$,  
the  
Lipschitz regularity    
and the uniform $L^\infty$ bound of $\Vatilde$,  
it follows from standard  
elliptic theories (see e.g., \cite{GT}, \cite{St} and  
\cite{Duff}) that  
$\Ga$ is uniquely defined by \eqref{Gaeq} and it satisfies  
(i) and (ii).\\  
Proof of (iii).     
Since $\Ga\in C_{\text{loc}}^2(M\setminus\{\xa\})$ we 
only need to check the inequality    
pointwise. If $\ua=0$ it is trivial. So assume $\ua>0$.    
By \eqref{Gaeq} we have     
\begin{align*}    
-\Delta_g\Ga+\Vatilde\Ga=0&&\text{pointwise in}\;M\setminus\{\xa\}.    
\end{align*}     
Since $\Ga>0$, using the definition of $\Vatilde$,    
we have    
\begin{equation*}    
\Delta_g\Ga=\Vatilde\Ga    
\le\big[c(n)\Rg+\alpha\big(\frac{\|\ua\|_{L^r(M)}}{\ua}\big)^{2-r}\big]\Ga,    
\end{equation*}    
pointwise in $(M\setminus\{\xa\})\cap\{\ua>0\}$.    
Multiplying the inequality above by $\ua$, we again obtain (iii).  
Proposition \ref{propGa} is established.
\end{proof}  
\begin{proof}[Proof of Proposition \ref{propuasupest}]  
The estimate for $\ua$ follows by \eqref{estimatenearxa} 
and Lemma \ref{lemzeaest}.  
Since $v_\alpha$ is uniformly bounded in $|y|<1$,
\eqref{uniformest} follows from the estimate of
 $\ua$ by a change of variables.
\end{proof}   

\section{Energy estimate}     
\label{secenergyest}     
We shall need estimates for the convergence rates of the limits ``in energy''
obtained in Section \ref{secprelims}.
The pointwise estimates obtained in Section 2 allow
us to adapt the energy estimates of  Bahri-Coron \cite{BC}.

 In order to simplify calculations, we
introduce a conformal metric $\hg=\psi^{4/(n-2)}g$,
with $\psi\in C^\infty(M)$,
$\psi(x_\alpha)=1$, $\frac 12 \le \psi\le 2$, $\|\psi\|_{ C^2}\le C$,
such that $R_{\hg}\equiv0$ in $B_{\delta_0}(P)$,
and where $\delta_0$ is a suitably chosen small constant and
both $\delta_0$ and $C$ depend only on $(M,g)$.
Such a metric may be obtained by locally solving
\[ 
-\Delg\psi+c(n)\Rg\psi=0
\qquad\text{in}\ B_{\delta_0}
\]
and then extending $\psi$ smoothly to $M$.
We denote, for $\delta_0/2\le \da\le \delta_0$,  
\begin{align*}
\tBa=B_{\da}^{\hg}(\xa).
\end{align*}

For $\tilde x\in \tBa $ and $\lambda>0$, we consider  
\[ 
\xi_{\tilde x,\lambda}^{\hg}(x)=\left(\frac{\lambda} 
{1+(\lambda\barl)^2\mathrm{dist}_{\hg}^2(x,\tilde x)^2}\right)^{\frac{n-2}{2}} 
\qquad\forall x\in\tBa. 
\]

It follows from Proposition \ref{propstandblowupva} that
\begin{equation}
\label{energyconv}
\lim_{\alpha\to+\infty}
\inttBa\{|\nabhg(\frac{\ua}{\psi}-\xi_{\xa,\mua^{-1}}^{\hg})|^2
+|\frac{\ua}{\psi}-\xi_{\xa,\mua^{-1}}^{\hg}|^{2^*}\}\dvhg=0.
\end{equation}

We follow the idea in \cite{BC} of selecting     
for every $\alpha$ an optimal multiple of a $\hg$-bubble, 
denoted $\ta\xia=\ta\xi_{\ya,\lda}^{\hg}$,  
and of estimating the difference $\ua/\psi-\ta\xia$
 by exploiting orthogonality.   
For future convenience, we prove our estimates for $n\ge3$.   
For $\tilde x\in B^{\hat g}_{\mua\da/2}(\xa)$ and $\lambda>0$,     
let $h_{\tilde x,\lambda}$ be defined by:     
\begin{equation}     
\label{haeq}     
\begin{cases}     
\Delhg h_{\tilde x,\lambda}=0\qquad&\text{in}\;\tBa\\     
h_{\tilde x,\lambda}=\xi_{\tilde x,\lambda}^{\hg}&\text{on}\;\de\tBa,     
\end{cases}     
\end{equation}     
and let $\chia$ be defined by     
\begin{equation}     
\label{chiaeq}     
\begin{cases}     
\Delhg\chia=0\qquad&\text{in}\;\tBa\\     
\chia=\frac{u_\alpha}\psi&\text{on}\;\de\tBa.     
\end{cases}     
\end{equation}     
Then $\ua/\psi-\chia\in H_0^1(\tBa)$, $\xi_{\tilde x,\lambda}^{\hg}-h_{\tilde x,
\lambda}
\in H_0^1(\tBa)$     
are the projections of $\ua$ and $\xi_{\tilde x,\lambda}^{\hg}$, respectively,  
on $H_0^1(\tBa)$.     
We set     
\[     
\sigma_{\tilde x,\lambda}=\xi_{\tilde x,\lambda}^{\hg}-h_{\tilde x,\lambda}. 
\]     
Then $\sigma_{\tilde x,\lambda}\le \xi_{\tilde x,\lambda}^{\hg}$ satisfies:     
\begin{equation*}     
\begin{cases}     
\Delhg\sigma_{\tilde x,\lambda}
=\Delhg\xi_{\tilde x,\lambda}^{\hg}\qquad&\text{in}\;\tBa\\     
\sigma_{\tilde x,\lambda}=0&\text{on}\;\de\tBa.     
\end{cases}     
\end{equation*}     
Let $(\ta,\ya,\lda)\in[\frac{1}{2},\frac{3}{2}]     
\times\overline{B^{\hat g}_{\mua\da/2}(\xa)}\times[\frac{1}{2\mua},\frac{3}{2\mua}]$      
be such that     
\begin{align*}     
\|\frac{\ua}{\psi}-\chia-&\ta\sigma_{\ya,\lda}\|_{\hg}\\    
=&\min     
\left\{\|\frac{\ua}{\psi}-\chia-t\sigma_{\tilde x,\lambda}\|_{\hg}:\;     
\begin{matrix}|t-1|\le1/2,\tilde x\in
\overline{B^{\hat g}_{\mua\da/2}(\xa)}\\     
|\mu_\alpha\lambda-1|\le1/2     
\end{matrix}     
\right\}.     
\end{align*}     
To simplify notation, henceforth we denote:    
\begin{align*}    
\sia=\sigma_{\ya,\lda},    
\qquad\xia=\xi_{\ya,\lda}^{\hg},    
\qquad\ha=h_{\ya,\lda},    
\end{align*}    
and we set:     
\[     
\wa=\frac{\ua}{\psi}-\chia-\ta\sia.     
\]   
The main result in this section is the following estimate 
for $\wa$:        
\begin{Proposition}[Energy estimate]   
\label{propwaestimate}     
For $n\ge3$, we have:     
\begin{align*}     
\|\wa\|+&|\ta^{2^*-2}\la-K^{-2}|\\   
\le&C(\mua^2\|U\|_{\stpr,\mua^{-1}}+\epsa\|U^{r-1}\|_{\stpr,\mua^{-1}} 
+\mua^{n-2}\|U^{2^*-2}\|_{\stpr,\mua^{-1}}).   
\end{align*}   
\end{Proposition}  
Recall from Section \ref{secprelims} that 
$\epsa=\mua^{n-\frac{n-2}{2}r}\alpha\|\ua\|_{L^r(M)}^{2-r}=\circ_\alpha(1)$. 
We define  
\begin{equation*}     
\Wa=\left\{w\in H_0^1(\tBa):     
\begin{matrix}     
&\bra\sia,w\ket_{\hg}=0\\     
&\bra f,w\ket_{\hg}=0   
\quad \forall f\in E   
\end{matrix}     
\right\},    
\end{equation*}    
where  
$E\subset  H_0^1(\tBa)$ is the tangent space
at $\sigma_{\ya,\lambda_\alpha}$ of   the finite dimensional surface  
$\{\sigma_{\tilde x,\lambda}:
 \tilde x\in B_{\mua\da}(\xa),\lambda>0\}\subset H_0^1(\tBa)$, 
with respect to the metric induced by  
the inner product $\bra u,v\ket_{\hg}=\inttBa\nabhg u\cdot\nabhg v\dvhg$. 
We work with coordinates given by the exponential 
map $\exp_{\ya}(y)$,
$y=(y^i)$,
 $i=1,\ldots,n$, we can write 
\[ 
E=\mathrm{span}\{\frac{\de\sia}{\de y^i},i=1,\ldots,n,\frac{\de\sia}{\de\lambda}\}, 
\] 
where 
\[ 
\frac{\de\sia}{\de y^i}=\frac{\de\sigma_{\exp_{\ya}(y),\lda}}{
\partial y^i}\big|_{y=0}, 
\qquad 
\frac{\de\sia}{\de\lambda}=\frac{\de\sigma_{\ya,\lambda}}{\de \lambda}\big|_{\lambda=\lda}. 
\] 
\begin{Lemma} 
\label{lemhanorm} 
For some constant $C$ independent of $\alpha$, 
\[ 
\inttBa|\nabhg\ha|^2\dvhg\le C\mua^{n-2}. 
\] 
\end{Lemma} 
\begin{proof} 
By standard elliptic estimates and properties of $\xia$,  
we have that 
$$
\inttBa|\nabhg\ha|^2\dvhg 
\le C\left(\int_{\de\tBa}|\nabhg\xia|^2\dshg+\int_{\de\tBa}\xia^2\dshg\right) 
\le C\mua^{n-2}. 
$$
\end{proof} 
We observe that by the uniform estimate \eqref{unifest} 
and by the maximum principle, 
\begin{equation} 
\label{maxhachia} 
\|\ha\|_{L^\infty(\tBa)}+\|\chia\|_{L^\infty(\tBa)}\le C\mua^{(n-2)/2}. 
\end{equation} 
It follows that $|\sia|\le C\xia$ on $\tBa$. 
Using Proposition \ref{propstandblowupva},     
it is not difficult to see that:     
\begin{Lemma}     
\label{lemtayala}     
As $\alpha\to+\infty$, we have   
$\|w_\alpha\|\to 0$, 
$\ta\to1$, $\mu_\alpha^{-1}
\mathrm{dist}_{\hat g}(x_\alpha, \ya)\to 0$, $\mu_\alpha\lambda_\alpha\to1$.     
Furthermore, $\wa\in\Wa$.     
\end{Lemma}     
\begin{proof}  
By definition of $\ta$ and $\sia$, 
\begin{align*} 
\|\ta\sia-\sigma_{\xa,\mua^{-1}}\| 
\le&\|\frac{\ua}{\psi}-\chia-\ta\sia\|+\|\frac{\ua}{\psi}-\chia-\sigma_{\xa,\mua^{-1}}\|\\ 
\le&2\|\frac{\ua}{\psi}-\chia-\sigma_{\xa,\mua^{-1}}\| 
\le2\|\frac{\ua}{\psi}-\xi^{\hat g}_{\xa,\mua^{-1}}\|. 
\end{align*} 
In the last step we have used  
$\Delhg(\frac{\ua}{\psi}-\xi^{\hat g}_{\xa,\mua^{-1}}) 
=\Delhg(\frac{\ua}{\psi}-\chia-\sigma_{\xa,\mua^{-1}})$ in $\tBa$ 
and $\frac{\ua}{\psi}-\chia-\sigma_{\xa,\mua^{-1}}=0$ on $\de\tBa$. 
Hence, in view of (\ref{energyconv}),
$\|\ta\sia-\sigma_{\xa,\mua^{-1}}\|\to0$ and
$\|w_\alpha\|\to 0$. 
By the arguments in Lemma \ref{lemhanorm}, 
we have 
$\|h_{\tilde x,\lambda}\|\le C\lambda^{-(n-2)/2}$ 
if $\mathrm{dist}_{\hg}(\tilde x,\xa)\le\da/2$. 
Consequently, we derive 
\begin{align*} 
\|\ta\xia-\xi^{\hat g}_{\xa,\mua^{-1}}\| 
\le\|\ta\sia-\sigma_{\xa,\mua^{-1}}\|+\|\ta\ha\|+\|h_{\xa,\mua^{-1}}\|\to0, 
\end{align*} 
as $\alpha\to+\infty$. 
It follows that
$\ta\to1$, $\mu_\alpha^{-1}
\mathrm{dist}_{\hat g}(x_\alpha, \ya)\to 0$, and $\mu_\alpha\lambda_\alpha\to1$.
Therefore the minimum of the norm is attained      
in the interior of      
$[\frac{1}{2},\frac{3}{2}]\times\overline{B_{\mua\da/2}(\xa)} 
\times[\frac{1}{2\mua},\frac{3}{2\mua}]$.     
Now a straightforward variational argument yields $\wa\in\Wa$.     
\end{proof}     
In order to estimate $\wa$, we begin by writing an equation     
for $\wa$:     
\begin{Lemma}     
\label{lemwaeq}     
$\wa$ satisfies:     
\begin{align}     
\label{wamaineq}     
-\Delhg\wa-\ka|\Thetaa|^{2^*-3} 
\Thetaa\wa+b'|\Thetaa|^{2^*-3}\wa^2+b''|\wa|^{2^*-1}=\fa 
\quad\text{in}\ \tBa, 
\end{align}  
where  
\begin{align*} 
&\ka=(2^*-1)\la\\ 
&\Thetaa=\ta\sia+\chia\\ 
&\fa=\la(\ta\xia)^{2^*-1}+\ta\Delhg\xia-\alpha\|\ua\|_{L^r(M)}^{2-r}\psi^{1-2^*}\ua^{r-1} 
+O(\mua^{(n-2)/2}\xia^{2^*-2}), 
\end{align*} 
and where $b',b''$ are bounded functions with $b'\equiv0$  
if $n\ge6$.      
\end{Lemma}     
\begin{proof}  
From \eqref{uaeq}, using the conformal invariance  
\eqref{conftransf} and recalling that $R_{\hg}\equiv0$ 
in $\tBa$, we have that $\ua/\psi$ satisfies: 
\[ 
-\Delhg\frac{\ua}{\psi} 
+\alpha\|\ua\|_{L^r(M)}^{2-r}\psi^{1-2^*}\ua^{r-1} 
=\la\big(\frac{\ua}{\psi}\big)^{2^*-1} 
\quad\text{in}\ \tBa. 
\] 
Consequently, $\wa$ satisfies: 
\begin{equation} 
\label{waprelimeq} 
-\Delhg\wa 
=\la(\Thetaa+\wa)^{2^*-1}+\ta\Delhg\sia 
-\alpha\|\ua\|_{L^r(M)}^{2-r}\psi^{1-2^*}\ua^{r-1} 
\qquad\text{in}\ \tBa. 
\end{equation} 
In order to simplify the right hand side in \eqref{waprelimeq},   
we use the elementary expansion:    
\begin{align*}     
(x+y)^{2^*-1}=|x|^{2^*-2}x&+(2^*-1)|x|^{2^*-3}xy\\   
&+b'(x,y)|x|^{2^*-3}y^2+b''(x,y)|y|^{2^*-1},     
\end{align*}     
for all $x,y\in\R$ such that $x+y\ge0$,     
where $b',b''$ are bounded functions and $b'\equiv0$  
if $n\ge6$.     
For $x=\Thetaa$ and $y=\wa$, we obtain:     
\begin{align*}    
(\Thetaa+\wa)^{2^*-1}=|\Thetaa|^{2^*-2}\Thetaa 
&+(2^*-1)|\Thetaa|^{2^*-3}\Thetaa\wa\\ 
&+b'|\Thetaa|^{2^*-3}\wa^2+b''|\wa|^{2^*-1}. 
\end{align*} 
Note that $\Thetaa=\ta\xia-\ta\ha+\chia$.       
By \eqref{maxhachia} and properties of $\xia$,  
we have $|\chia-\ta\ha|\le C\mua^{(n-2)/2} 
\le C\ta\xia$, and thus     
by simple calculus:     
\begin{align*}     
|\Thetaa|^{2^*-2}\Thetaa=&|\ta\xia-\ta\ha+\chia|^{2^*-2}(\ta\xia-\ta\ha+\chia)\\ 
&=(\ta\xia)^{2^*-1}+O(\mua^{(n-2)/2}\xia^{2^*-2}). 
\end{align*}      
Inserting the above expansions into \eqref{waprelimeq},  
we obtain     
\eqref{wamaineq}.     
\end{proof}     
The proof of Proposition \ref{propwaestimate}     
relies on the coercivity property as in Lemma \ref{lemcoercivity} below.    
Recall that $\delta_0>0$ was introduced in Section \ref{secprelims}    
as an upper bound for the radii $\da$ of the balls $\tBa=B_{\da}(\xa)$.     
Here is where we fix $\delta_0$.    
We denote by $\Qa$ the continuous bilinear form defined 
for $\varphi,\psi\in H_0^1(\tBa)$ by: 
\[ 
\Qa(\varphi,\psi)=\inttBa\{\nabhg\varphi\cdot\nabhg\psi 
-\ka|\Thetaa|^{2^*-3}\Thetaa\,\varphi\psi\}\dvhg, 
\] 
where $\ka$ and $\Thetaa$ are defined in Lemma \ref{lemwaeq}. 
\begin{Lemma}     
\label{lemcoercivity}   
There exist $0<\delta_0\ll1$, $\alpha_0\gg1$    
and $c_0>0$ independent of $\alpha\ge\alpha_0$ such that    
\begin{align*}     
\Qa(w,w)\ge c_0\inttBa|\nabhg w|^2\dvhg,    
\qquad\forall w\in\Wa,\ \forall\alpha\ge\alpha_0.    
\end{align*}    
\end{Lemma}    
Lemma \ref{lemcoercivity} is a consequence of the following 
general perturbation result: 
\begin{Lemma} 
\label{lemperturbation} 
Let $\Omega\subset\R^n$, let $h$ be a metric on $\Omega$, 
$k>0$ and $\Theta\in L^{2^*}(\Omega)$.  
Denote by $Q$ the continuous bilinear form defined on 
$H_0^1(\Omega)\times H_0^1(\Omega)$ by 
\[ 
Q(\varphi,\psi)=\int_{\Omega}\{\nabla_h\varphi\cdot\nabla_h\psi 
-k|\Theta|^{2^*-3}\Theta\,\varphi\psi\}\,dv_h. 
\] 
There exist $\eps_0>0$ and $c_1>0$, depending only 
on $n$,  such that if 
\[ 
\|\Theta-U\|_{L^{2^*}(\Omega)}+|k-(2^*-1)K^{-2}|+
\|h-\mathcal E\|_{L^\infty(\Omega)}
\le\eps_0, 
\] 
where $\mathcal E$ denotes the Euclidean metric, 
then 
\[ 
Q(\varphi,\varphi)\ge\frac{c_1}{2}\int_\Omega|\nabla_h\varphi|^2\,dv_h, 
\qquad\begin{matrix} 
&\forall\varphi\in H_0^1(\Omega):\\ 
&|\bra\varphi,e_i\ket_h|\le\eps_0\|\varphi\|_h,\ 0\le i\le n+1, 
\end{matrix} 
\] 
where $e_0=U$, $e_i=\de U_{y,1}/\de y^i\big|_{y=0}$, $i=1,\ldots,n$,  
$e_{n+1}=\de U_{0,\lambda}/\de\lambda\big|_{\lambda=1}$.  
\end{Lemma} 
\begin{proof} 
We denote by $\widetilde Q$ the continuous bilinear form on $D^{1,2}(\R^n)$    
defined by    
\begin{align*}    
\widetilde Q(\widetilde\varphi,\widetilde\psi)    
=\intRn\{\nabla\widetilde\varphi\cdot\nabla\widetilde\psi    
-\frac{2^*-1}{K^2}U^{2^*-2}\,\widetilde\varphi\widetilde\psi\}\dy.    
\end{align*}    
It is well-known (and it may be verified    
by pull-back to the standard $n$-sphere,  
in stereographic projection coordinates)    
that there exists $c_1>0$ such that    
\begin{align}    
\label{flatcoercivity}    
\widetilde Q(\widetilde\varphi,\widetilde\varphi)    
\ge c_1\intRn|\nabla\widetilde\varphi|^2\dy,    
\qquad\begin{matrix}\forall\widetilde\varphi\in D^{1,2}(\R^n):\\  
\langle\widetilde\varphi,e_i\rangle_{\mathcal E}=0,\    
0\le i\le n+1.    
\end{matrix}  
\end{align} 
Now the claim follows by elementary considerations. Indeed, there exist
unique $\mu^j$, $|\mu^j|=O(\eps_0\|\varphi\|_h)$, such that
$\widetilde\varphi:=\varphi -\mu^j e_j$ satisfies
$\langle\widetilde\varphi, e_i\rangle_{\mathcal E}=0$, 
$\forall\ 0\le i\le n+1$.  (\ref{flatcoercivity}) holds for
$\widetilde\varphi$, and the claim follows easily. 
\end{proof}

We introduce some notations:   We set
\[
\tOa=\mua^{-1}(\exp_{\ya}^{\hg})^{-1}(\tBa)\subset\R^n.
\]
We denote by $\Tmua$ the transformation
which maps $f:\tBa\to\R$ into $\Tmua f:\tOa\to\R$ defined by
\[
(\Tmua f)(y)=\mua^{(n-2)/2}f(\exp_{\ya}^{\hg}(\mua y))
\qquad\forall y\in\tOa.
\]
We denote by $\ga$ the metric on $\tOa$ defined by
$\ga(y)=\hg(\exp_{\ya}^{\hg}(\mua y))$.
The following transformation properties hold:
\begin{align}
\label{gradtransf}
\inttBa\nabhg\varphi\cdot\nabhg  \psi \dvhg
=\intta\nabga\Tmua\varphi\cdot\nabga\Tmua
\psi \dvga
\qquad\forall\varphi,  \psi \in H_0^1(\tBa).
\end{align}
If $p_1+\cdots+p_k=2^*$, then
\begin{align}
\label{prodtransf}
\inttBa|\varphi_1|^{p_1}\cdots|\varphi_k|^{p_k}\dvhg
=&\intta|\Tmua\varphi_1|^{p_1}\cdots|\Tmua\varphi_k|^{p_k}\dvga,\\
\nonumber
&\forall \varphi_1,\ldots,\varphi_k\in H_0^1(\tBa).
\end{align}

\begin{proof}[Proof of Lemma \ref{lemcoercivity}]   
Observe that by the transformation properties 
\eqref{gradtransf}--\eqref{prodtransf} we have 
\[ 
\Qa(\varphi,\varphi) 
=\intta\{|\nabga\Tmua\varphi|^2 
-\ka|\Tmua\Theta|^{2^*-3}\Tmua\Theta\,(\Tmua\varphi)^2\}\dvga. 
\] 
By taking $\delta_0$ small, we achieve $|\ga-\mathcal E|\le\eps_0$.  
By taking $\alpha_0\gg1$,  
we achieve $\|\Tmua\Thetaa-U\|_{L^{2^*}(\tOa)}\le\eps_0$ and $|\ka-(2^*-1)K^{-2}|\le\eps_0$. 
It remains to check that by taking a possibly smaller $\delta_0$ and a possibly larger $\alpha_0$, 
we have for all $\varphi\in\Wa$: 
\begin{align*} 
|\bra\Tmua\varphi,e_i\ket_{\mathcal E}|\le\eps_0\|\Tmua\varphi\|_{\ga}. 
\end{align*} 
We check the above for $i=0$. 
Since $\varphi\in\Wa$, 
\begin{align*} 
0=&\inttBa\nabhg\varphi\cdot\nabhg\sia\dvhg 
=\inttBa\nabhg\varphi\cdot\nabhg\xia\dvhg\\ 
=&\intta\nabga\Tmua\varphi\cdot\nabga\Tmua\xia\dvga 
=\intta\ga^{ij}\frac{\de\Tmua\varphi}{\de y^i} 
\frac{\de\xia}{\de y^j}\sqrt{\det\ga}\dy. 
\end{align*} 
Therefore: 
\begin{align*} 
|&\intta\nabla\Tmua\varphi\cdot\nabla U\dy|\\ 
\le&|\intta\nabla\Tmua\varphi\cdot\nabla(U-\Tmua\xia)\dy| 
+|\intta(\delta^{ij}-\ga^{ij}\sqrt{\det\ga}) 
\frac{\de\Tmua\varphi}{\de y^i}\frac{\de\xia}{\de y^j}\dy|\\ 
\le&(\|U-\Tmua\xia\|+\sup_{\tOa}|\delta^{ij}-\ga^{ij}\sqrt{\det\ga}|\|\nabla\Tmua\xia\|) 
\times\|\Tmua\varphi\|\\ 
=&(\circ_\alpha(1)+O_{\delta_0}(\delta_0^2))\|\Tmua\varphi\|_{\ga}. 
\end{align*} 
The remaining conditions are verified similarly. 
Taking into account \eqref{gradtransf}, we conclude 
by Lemma \ref{lemperturbation} that 
for all $\varphi\in\Wa$, 
\[ 
\Qa(\varphi,\varphi)\ge\frac{c_1}{2}\intta|\nabga\Tmua\varphi|^2\dvga 
=\frac{c_1}{2}\inttBa|\nabhg\varphi|^2\dvhg, 
\] 
as asserted. 
\end{proof}    
Now the     
\begin{proof}[Proof of Proposition \ref{propwaestimate}]   
Multiplying \eqref{wamaineq} by $\wa$ and integrating over $\tBa$     
we obtain:     
\begin{align*}     
\Qa(\wa,\wa)+\circ_\alpha(\|\wa\|^2) 
=\inttBa\fa\wa\dvhg. 
\end{align*} 
By Lemma \ref{lemcoercivity},    
in view of the form of $\fa$ and recalling the orthogonality property 
$\inttBa\wa\Delhg\xia\dvhg=0$, we derive from the above: 
\begin{align*} 
\|\wa\|\le C(\|\xia^{2^*-1}&+K^2\Delhg\xia\|_{\stpr}\\ 
&+\alpha\|\ua\|_{L^r(M)}^{2-r}\|\ua^{r-1}\|_{\stpr} 
+\mua^{(n-2)/2}\|\xia^{2^*-2}\|_{\stpr}). 
\end{align*} 
$\xia$ satisfies 
\begin{equation} 
\label{xiaeq} 
-\Delhg\xia=K^{-2}\xia^{2^*-1}+O(\xia). 
\end{equation} 
It follows that 
\[ 
\|\xia^{2^*-1}+K^2\Delhg\xia\|_{\stpr} 
\le C\|\xia\|_{\stpr} 
\le C\mua^2\|U\|_{\stpr,\mua^{-1}}, 
\] 
where we have used $\lda\mua\le C$ in the last inequality. 
In order to estimate the second term, we note that 
the uniform estimate \eqref{unifest} implies: 
$\ua\le C   \xia  $ in $\tBa$. 
Consequently, 
\begin{align*} 
\|\ua^{r-1}\|_{\stpr} 
\le C\left(\inttBa\xi_{\ya,\mua^{-1}}^{(r-1)\stpr}\dvhg\right)^{1/\stpr} 
\le C\mua^{n-\frac{n-2}{2}r}\|U^{r-1}\|_{\stpr,\mua^{-1}}. 
\end{align*} 
Similarly, we compute: 
\[ 
\|\xia^{2^*-2}\|_{\stpr} 
\le C\mua^{(n-2)/2}\|U^{2^*-2}\|_{\stpr,\mua^{-1}}. 
\] 
The asserted decay estimate for $\|\wa\|$ follows.\par  
In order to estimate $|\la\ta^{2^*-1}-\ta K^{-2}|$, 
in view of \eqref{xiaeq}, 
we write $\fa$ in the form: 
\begin{align*} 
\fa=&(\la\ta^{2^*-1}-\ta K^{-2})\xia^{2^*-1} 
-\alpha\|\ua\|_{L^r(M)}^{2-r}\psi^{1-2^*}\ua^{r-1}\\ 
&+O(\xia)+O(\mua^{(n-2)/2}\xia^{2^*-2}). 
\end{align*} 
Multiplying \eqref{wamaineq} by $\sia$, integrating over $\tBa$   
and taking into account that $\inttBa\sia\Delhg\wa\dvhg=0$,  
we have:  
\begin{align*} 
&-\ka\inttBa|\Thetaa|^{2^*-3}\Thetaa\wa\sia\dvhg 
+\inttBa b'|\Theta|^{2^*-3}\wa^2\sia\dvhg\\ 
&\qquad+\inttBa b''|\wa|^{2^*-1}\sia\dvhg\\ 
&=(\la\ta^{2^*-1}-\ta K^{-2})\inttBa\xia^{2^*-1}\sia\dvhg 
-\alpha\|\ua\|_{L^r(M)}^{2-r}\inttBa\psi^{1-2^*}\ua^{r-1}\sia\dvhg\\ 
&\qquad+O(\inttBa\xia\sia\dvhg)+O(\mua^{(n-2)/2})\inttBa\xia^{2^*-2}\sia\dvhg 
\end{align*} 
and thus, using $|\wa|+|\Thetaa|
+|\sia|+\ua\le C\xia$, we derive: 
\begin{align*} 
|\la\ta^{2^*-1}&-\ta K^{-2}|\inttBa\xia^{2^*-1}\sia\dvhg\\ 
\le&C\big(\inttBa\xia^{2^*-1}|\wa|\dvhg 
+\inttBa\xia^2\dvhg\\ 
&+\alpha\|\ua\|_{L^r(M)}^{2-r}\inttBa\xia^{r}\dvhg 
+\mua^{(n-2)/2}\inttBa\xia^{2^*-1}\dvhg\big). 
\end{align*} 
In order to compare with the decay rate of $\|\wa\|$, it is 
convenient to estimate as follows: 
\[ 
\inttBa\xia^{2^*-1}|\wa|\dvhg\le C\|\wa\|; 
\] 
$$
\inttBa\xia^{r}\dvhg 
\le C
\|\xia^{r-1}\|_{\stpr}\|\xia\|_{2^*}
\le C\mua^{n-\frac{n-2}{2}r}\|U^{r-1}\|_{\stpr,\mua^{-1}};
$$
$$
\mua^{(n-2)/2}\inttBa\xia^{2^*-1}\dvhg 
=\mua^{(n-2)/2}\|\xia^{2^*-2}\|_{\stpr}\|\xia\|_{2^*}
\le C\mua^{n-2}\|U^{2^*-2}\|_{\stpr,\mua^{-1}};
$$
\[ 
\inttBa\xia^2\dvhg 
\le C\|\xia\|_{2^*}\|\xia\|_{\stpr}\le C\mua^2\|U\|_{\stpr,\mua^{-1}}. 
\] 
On the other hand, 
\begin{align*} 
|\la\ta^{2^*-1}&-\ta K^{-2}|\inttBa\xia^{2^*-1}\sia\dvhg\\ 
\ge&|\la\ta^{2^*-1}-\ta K^{-2}|\big(\inttBa\xia^{2^*}\dvhg+O(\mua^{n-2})\big)\\ 
\ge&C^{-1}|\la\ta^{2^*-1}-\ta K^{-2}|+O(\mua^{n-2}). 
\end{align*} 
The estimate for $|\la\ta^{2^*-1}-\ta K^{-2}|$ is established.    
\end{proof}    
\section{Lower bound for $\Yg$ and proof of Theorem \ref{thmmain}  
for $n\ge7$}    
\label{secsplitting}   
In this section we shall carefully exploit orthogonality  
in order to derive a lower bound for $\Yg(\ua)$,  
as in Proposition \ref{propyamabelowerbound} below.  
Together with the estimates from the previous sections,  
it will readily imply the proof of Theorem \ref{thmmain}  
in the case $n\ge7$.   
We shall need an $L^2$-estimate of $|\nabg\ua|$ on $\de\tBa$.    
This can be achieved by selecting a suitable ``good radius"   
$\da\in[\delta_0/2,\delta_0]$, see Lemma \ref{lemgradbdry} below.    
Here is where we fix $\da$.  
Unless otherwise stated, we assume $n\ge3$.\par  
The main step towards obtaining a contradiction is the following     
lower bound for $\Yg(\ua)$:    
\begin{Proposition}[Lower bound for $\Yg$]    
\label{propyamabelowerbound}    
Let $\delta_\alpha$ be a ``good radius".    
Then, for all $\alpha$ sufficiently large,  
\begin{align*} 
\Yg(\ua)\ge Y_{\hg}(\widetilde\xi_{\ya,\lda}^{\hg}) 
+O(\mua^2\|U\|_{\stpr,\mua^{-1}}\|\wa\|+\mua^{n-2}). 
\end{align*}  
\end{Proposition}   
Proposition \ref{propyamabelowerbound} readily implies: 
\begin{Corollary} 
\label{corlowerbound} 
The following estimates hold: 
\begin{align*} 
\tag{i} 
\alpha\|\ua\|_{L^r(M)}^2\le&C(\mua^2\|U\|_{\stpr,\mua^{-1}}+\mua^{n-2}) 
+|K^{-2}-Y_{\hg}(\widetilde\xi_{\ya,\lda}^{\hg})|\\ 
\tag{ii} 
\alpha\|\ua\|_{L^r(M)}^2 
\le&C(\mua^4\|U\|_{\stpr,\mua^{-1}}^2 
+\epsa\mua^2\|U\|_{\stpr,\mua^{-1}}\|U^{r-1}\|_{\stpr,\mua^{-1}}\\ 
&\quad+\mua^n\|U\|_{\stpr,\mua^{-1}}\|U^{2^*-2}\|_{\stpr,\mua^{-1}}+\mua^{n-2})\\ 
&\quad+|K^{-2}-Y_{\hg}(\widetilde\xi_{\ya,\lda}^{\hg})|. 
\end{align*} 
\end{Corollary}  
\begin{proof} 
By the initial assumption on $\Ia$, we have: 
\[ 
K^{-2}>\Ia(\ua)=\Yg(\ua)+\alpha\|\ua\|_{L^r(M)}^2. 
\] 
Therefore, the lower bound as in Proposition \ref{propyamabelowerbound} implies: 
\[ 
K^{-2}>\alpha\|\ua\|_{L^r(M)}^2+Y_{\hg}(\widetilde\xi_{\ya,\lda}^{\hg}) 
+O(\mua^2\|U\|_{\stpr,\mua^{-1}}\|\wa\|+\mua^{n-2}), 
\] 
which in turn yields: 
\[ 
\alpha\|\ua\|_{L^r(M)}^2\le|Y_{\hg}(\widetilde\xi_{\ya,\lda}^{\hg})-K^{-2}| 
+C(\mua^2\|U\|_{\stpr,\mua^{-1}}\|\wa\|+\mua^{n-2}). 
\] 
Using $\|\wa\|\le C$, we obtain (i). 
Using the energy estimate as in Proposition \ref{propwaestimate}, 
we obtain (ii). 
\end{proof} 
The proof of Proposition \ref{propyamabelowerbound}    
relies on some boundary estimates and on consequences of orthogonality,    
which we proceed to derive.    
\begin{Lemma}[Choice of ``good radius"]    
\label{lemgradbdry}    
There exists $C>0$ independent of $\alpha$    
such that:    
\[    
\int_{M\setminus B_{\delta_0/2}(\xa)}|\nabg\ua|^2\dvg\le C\mua^{n-2}.    
\]    
Consequently, for every $\alpha$ we can select $\da\in[\delta_0/2,\delta_0]$    
such that on $\tBa=B_{\da}(\xa)$ we have:    
\begin{align}   
\label{goodradiusest}   
\int_{\de\tBa}|\nabg\ua|^2\dsg\le C\mua^{n-2}.  
\end{align} 
Furthermore, for such a $\da$ we have: 
\begin{equation} 
\label{chiabdry} 
\inttBa|\nabhg\chia|^2\dvhg\le C\mua^{n-2}. 
\end{equation}    
\end{Lemma}    
\begin{proof}    
Denote by $\eta$ a smooth cutoff function to be fixed below, 
satisfying $0\le\eta\le1$.    
Multiplying \eqref{uaeq} by $\eta^2\ua$ and integrating    
by parts on $M$ we have:    
\[    
\intM\nabg\ua\cdot\nabg(\eta^2\ua)\dvg\le-c(n)\intM\Rg\eta^2\ua^2\dvg    
+\la\intM\eta^2\ua^{2^*}\dvg.    
\]    
It follows that:    
\[    
\intM\eta^2|\nabg\ua|^2\dvg\le C\big(\intM\ua^2(|\nabg\eta|^2+\eta^2)\dvg    
+\intM\eta^2\ua^{2^*}\dvg\big).    
\]    
Choosing $\eta$ such that $\eta\equiv1$ in     
$M\setminus B_{\delta_0}(\xa)$,    
$\text{supp}\,\eta\subset M\setminus B_{\delta_0/2}(\xa)$,    
we obtain:    
\begin{align*}    
\int_{M\setminus B_{\delta_0}(\xa)}|\nabg\ua|^2\dvg    
\le C\big(\int_{M\setminus B_{\delta_0/2}(\xa)}\ua^2\dvg    
+\int_{M\setminus B_{\delta_0/2}(\xa)}\ua^{2^*}\dvg\big).    
\end{align*}    
Now the statement follows    
by the uniform estimate as in Proposition \ref{propuasupest}.\\    
\noindent    
Since $\ua\in C^1(M)$, we can choose $\da$ such that:    
\[    
\int_{\de\tBa}|\nabg\ua|^2\dvg    
=\min_{\delta\in[\delta_0/2,\delta_0]}  
\int_{\de B_\delta(\xa)}|\nabg\ua|^2\dsg\le C\mua^{n-2}.    
\]   
Recalling the definition of $\chia$, we have by standard 
elliptic estimates and equivalence of $g$ and $\hg$: 
\[ 
\int_{\de\tBa}|\nabhg\chia|^2\dvhg 
\le C\int_{\de\tBa}\{|\nabg \ua|^2+\ua^2\}\dsg 
\le C\mua^{n-2}. 
\]  
\end{proof} 
\begin{proof}[Proof of Proposition \ref{propyamabelowerbound}]  
By the uniform estimate \eqref{unifest} and by
Lemma \ref{lemgradbdry},  
we have 
\begin{equation} 
\label{dropbdry} 
\Yg(\ua)=\frac{\inttBa\{|\nabg\ua|^2+c(n)\Rg \ua^2\}\dvg} 
{\big(\inttBa\ua^{2^*}\dvg\big)^{2/2^*}}+O(\mua^{n-2}). 
\end{equation} 
By conformal invariance \eqref{conftransf},
together with \eqref{goodradiusest} and \eqref{unifest},  
\[ 
\inttBa\{|\nabg\ua|^2+c(n)\Rg \ua^2\}\dvg 
=\inttBa|\nabhg\frac{\ua}{\psi}|^2\dvhg+O(\mua^{n-2}). 
\] 
Recall from Section \ref{secenergyest} 
that $\ua/\psi=\ta\xia-\ta\ha+\chia+\wa$. 
By Lemma \ref{lemhanorm}, \eqref{maxhachia}, \eqref{chiabdry},
and the fact
\[
\intBa\nabhg\ha\cdot\nabhg\wa\dvhg
=0=\intBa\nabhg\chia\cdot\nabhg\wa\dvhg,
\]
we have
\begin{equation}
\label{drophachia}
\Yg(\ua)= F(\wa)+O(\mua^{n-2}),
\end{equation}
where
\[
F(w):=\frac{\intBa|\nabhg (\ta\xia+w)|^2\dvhg}
{\big(\intBa|\ta\xia+w|^{2^*}\dvhg\big)^{2/2^*}},
\qquad w\in H^1_0(\Ba).
\]
A Taylor expansion yields:
\begin{align}
\label{taylor}
F(\wa)=F(0)+F'(0)\wa+\frac{1}{2}\bra F''(0)\wa,\wa\ket
+\circ(\|\wa\|^2),
\end{align}
where $F'$, $F''$ denote Fr\'echet derivatives.
We compute:
\begin{align*}
F'&(0)\wa=
\frac{2}{(\intBa(\ta\xia\dvhg)^{2^*})^{2/2^*}}\times\\
&\times\big\{\intBa\nabhg(\ta\xia)\cdot\nabhg\wa\dvhg
-\frac{\intBa|\nabhg\ta\xia|^2\dvhg}{\intBa(\ta\xia)^{2^*}\dvhg}
\intBa(\ta\xia)^{2^*-1}\wa\dvhg\big\}.
\end{align*}
By orthogonality,
$\intBa\nabhg\xia\cdot\nabhg\wa\dvhg=0$
and by \eqref{xiaeq}
\[
K^2\intBa\xia^{2^*-1}\wa\dvhg
=\intBa(-\Delhg\xia+O(\xia))\wa\dvhg
=O(\intBa\xia\wa\dvhg).
\]
Hence,
\[
|\intBa\xia^{2^*-1}\wa\dvhg|
\le C\|\xia\|_{\stpr}\|\wa\|
\le C\mua^2\|U\|_{\stpr,\mua^{-1}}\|\wa\|
\]
and consequently
\begin{equation*}
|F'(0)\wa|\le C\mua^2\|U\|_{\stpr,\mua^{-1}}\|\wa\|.
\end{equation*}
Similarly, we compute:
\begin{align*}
\bra F''&(0)\wa,\wa\ket
=\frac{2}{\big(\intBa(\ta\xia)^{2^*}\big)^{2/2^*}}\times\\
&\times\big\{\intBa|\nabhg\wa|^2\dvhg
-(2^*-1)\frac{\intBa|\nabhg\xia|^2\dvhg}{\intBa\xia^{2^*}\dvhg}
\intBa\xia^{2^*-2}\wa^2\dvhg\big\}\\
&\qquad+O(\intBa\xia^{2^*-1}\wa\dvhg)^2.
\end{align*}
By the transformations
\eqref{gradtransf}--\eqref{prodtransf} and by Lemma \ref{lemperturbation}
with
\begin{align*}
&\Omega=\tOa,\qquad\Theta=\Tmua\xia,\qquad h=\ga\\
&k=(2^*-1)\frac{\inttBa|\nabhg\xia|^2\dvhg}{\inttBa\xia^{2^*}\dvhg},
\end{align*}
we obtain, for large $\alpha$, that
\begin{align*}
\intBa|\nabhg\wa|^2\dvhg
-(2^*-1)\frac{\intBa|\nabhg\xia|^2\dvhg}{\intBa\xia^{2^*}\dvhg}
\intBa\xia^{2^*-2}\wa^2\dvhg\ge\frac{c_1}{2}\|\wa\|^2.
\end{align*}
Consequently,
\begin{align*}
\bra F''(0)\wa,\wa\ket
\ge\frac{c_1}{2}\|\wa\|^2+O(\mua^4)\|U\|_{\stpr,\mua^{-1}}^2\|\wa\|^2.
\end{align*}
Inserting into \eqref{taylor} and observing that
$\mua^2\|U\|_{\stpr,\mua^{-1}}=\circ_\alpha(1)$, we derive:
\begin{equation}
\label{lowerbd}
F(\wa)\ge F(0)
+O(\mua^2\|U\|_{\stpr,\mua^{-1}}\|\wa\|+\mua^{n-2}).
\end{equation}
Returning to \eqref{dropbdry}
and taking into account that
\[
F(0)=Y_{\hg}(\tilde\xi_{\tilde\xa,\lda}^{\hg})
+O(\mua^{n-2}),
\]
we obtain the asserted lower bound.
\end{proof}
\begin{proof}[Proof of Theorem \ref{thmmain} for $n\ge7$]  
By straightforward computations, 
\[ 
\|U^q\|_{\stpr,\mua^{-1}}\le 
\begin{cases} 
C,&\text{if}\ q>(n+2)/[2(n+2)]\\ 
(\log\mua^{-1})^{1/\stpr} 
&\text{if}\ q=(n+2)/[2(n+2)]\\ 
\mua^{-\frac{n+2}{2}+q(n-2)} 
&\text{if}\ q<(n+2)/[2(n+2)]\\ 
\end{cases}. 
\] 
We take $r=\br=2n/(n+2)$.  
Then, since $n\ge7$, we have:  
\begin{align*}  
&\|U\|_{\stpr,\mua^{-1}}\le C\\  
&\|U^{\br-1}\|_{\stpr,\mua^{-1}}     
\le C(1+\mua^{-2+\beta})\\   
&\|U^{2^*-2}\|_{\stpr,\mua^{-1}}\le C\mua^{-(n-6)/2},   
\end{align*}  
where $\beta=(n-6)(n-2)/[2(n+2)]$ is {\em strictly positive}.  
Hence, (ii) in Corollary \ref{corlowerbound} yields:  
\begin{equation}  
\label{7alphaupperbound}  
\alpha\|\ua\|_{L^{\br}(M)}^2\le|Y_{\hg}(\xia,\tBa)-K^{-2}|  
+C[\mua^4+\epsa(\mua^2+\mua^\beta)].  
\end{equation}  
By \eqref{aubinexpansion},   
\[  
|Y_{\hg}(\widetilde\xi_{\ya,\lda}^{\hg})-K^{-2}|\le C\mua^4.  
\]  
In view of \eqref{epsa}, we derive:  
\begin{equation}  
\label{alphaupperbound} 
\alpha\|\ua\|_{L^{\br}(M)}^2\le C\mua^4. 
\end{equation}  
On the other hand, rescaling, we have:  
\begin{equation}  
\|\ua\|_{L^{\br}(M)}\ge\|\ua\|_{L^{\br}(\tBa)}  
\ge C^{-1}\mua^2\|U\|_{L^{\br}(B_1(0))}\ge C^{-1}\mua^2  
\end{equation}  
and inserting into \eqref{alphaupperbound} we obtain   
$\alpha\le C$, a contradiction.    
Hence, Theorem \ref{thmmain} is established for all $n\ge7$.   
\end{proof}  
\section{Proof of Theorem \ref{thmmain}
for  $n=6$}  
\label{seclimitcase}  
In order to prove Theorem \ref{thmmain}  
in the remaining case $n=6$ we   
need a uniform {\em lower} bound for $\ua$.  
Indeed we shall prove: 
\begin{Proposition}[Uniform lower bound]   
\label{proplowerbound}    
For $n=6$, $r=\br=3/2$, and any 
$1/2<\gamma<1$,  there exists some   constant
 $C>0$,  which is independent of $\alpha$, such that $\ua$ satisfies:    
\[    
\ua(x)\ge C^{-1}\mua^2\distg(x,\ya)^{-4}    
\qquad\forall x\in B_{\da\mua^{\gamma}}(\xa)\setminus B_{\mua}(\xa),    
\]    
for all $\alpha\gg1$.    
\end{Proposition}    
\begin{proof}    
We equivalently show that 
\[ 
\va(y)\ge C^{-1}|y|^{-4}, 
\qquad\forall y\in B_{\da\mua^{\gamma-1}}(0)\setminus B_{1}. 
\] 
Here $\va$ is defined on $\Oa$ as in (\ref{100}) and
(\ref{200}).   Recall the $\delta_0/2\le \delta_\alpha\le \delta_0$.
$0<\delta_0<1$ will be small and fixed below.
We define a comparison function  
\[    
\Ha(y)=\tau\,\big(\frac{\da^4}{|y|^4}-\mua^4\big)    
+L\,\mua^2\big( \log \frac 1{\mu_\alpha}\big)^{2/3}
\log\big(\frac{\mua|y|}{\da}\big),   
\qquad y\in\Oa\setminus B_1  
\]    
where $\tau>0$, $L>0$ will be chosen below.

Since $\va\to U$ uniformly on $\partial B_1$, 
we first fix some $0<\tau=\tau(\delta_0)<1$ such that
$\va\ge H_\alpha$ on $\partial B_1$ for large $\alpha$.
Since $H_\alpha=0$ on $\partial \Oa$, we also have
$\va\ge H_\alpha$ on $\partial \Oa$.
We know that
$$
 C_1^{-1}\delta_0^{-2}\mua^2\le |y|^{-2}\le 1,\qquad \mbox{on}\
\Oa\setminus B_1.
$$
Here and in the following, $ C_1>1$ denotes some constant
depending only on $(M,g)$.
Setting $\ga(y)=g(\exp_{\xa}^{g}(\mua y))$ 
we have, 
\begin{align*} 
&|\Delga|y|^{-4}|\le  C_1\mua^2|y|^{-4}\\ 
&|\Delga\log|y|-4|y|^{-2}|\le  C_1\mua^2. 
\end{align*} 
Hence,
$$
\Delga H_\alpha(y)
\ge 4L\mua^2(\log\frac{1}{\mua})^{2/3}|y|^{-2}
- C_1 \mua^2|y|^{-4}
- C_1 L \mua^4(\log\frac{1}{\mua})^{2/3}.
$$
Recall from Section 1 that 
$$
\eps_\alpha
=\alpha\mua^{n-\frac{n-2}{2}r}\|\ua\|_{L^r(M)}^{2-r}
\le \alpha\|\ua\|_{L^r(M)}^2.
$$
By (i) in Corollary \ref{corlowerbound},
$$
\alpha\|\ua\|_{L^r(M)}^2\le C_2
(\mua^2\|U\|_{\stpr,\mua^{-1}}+\mua^4)
+|Y_{\hg}(\widetilde\xi_{\ya,\lda}^{\hg})-K^{-2}|.
$$
Here and in the following, $C_2>1$ denotes some constant independent
of $\alpha$ and $L$.
By the expansion \eqref{aubinexpansion},
\begin{equation}
|Y_{\hg}(\widetilde\xi_{\ya,\lda}^{\hg})-K^{-2}|\le C_2\mua^4\log\frac1\mua,
\label{300}
\end{equation} 
and, clearly, 
$$
\|U\|_{\stpr,\mua^{-1}}\le C_2(\log\frac{1}{\mua})^{2/3}.
$$
It follows that
$$
\alpha\|\ua\|_{L^r(M)}^2\le C_2
\mua^2\big( \log\frac1\mua \big)^{2/3}.
$$
Together with the uniform estimate from 
Section 2: $v_\alpha(y)\le C_2|y|^{-4}$, we obtain
$$
\eps_\alpha \va^{1/2}\le C_2 \mua^2 \big( \log\frac1\mua \big)^{2/3}
|y|^{-2}.
$$
Hence, by the equation of $\va$,
we have, on $\Oa\setminus B_1$, that
\begin{align*} 
&-\Delga(\va-\Ha)(y)\\
\ge& (4L-C_2)\mua^2 (\log\frac{1}{\mua})^{2/3}|y|^{-2}
-C_2\mua^2|y|^{-4}-C_1 L \mua^4 (\log\frac{1}{\mua})^{2/3}\\
\ge& (4L-C_2-C_1\delta_0^2L)\mua^2 (\log\frac{1}{\mua})^{2/3}|y|^{-2}.
\end{align*} 
We first fix $\delta_0>0$ small
($C_1\delta_0<1$), and then
 take $L$ large,  we achieve, for large $\alpha$, that 
\[ 
-\Delga(\va-\Ha)\ge0 
\qquad\text{in}\ \Oa\setminus B_1. 
\] 
By the maximum principle, 
\[ 
\va\ge\Ha\qquad\text{in}\ \Oa\setminus B_{1}. 
\] 
To conclude, we observe that for any fixed $1/2<\gamma<1$ 
we can find a $C>0$ such that: 
\[ 
\Ha(y)\ge C^{-1}|y|^{-4} 
\qquad\text{in}\ B_{\da\mua^{\gamma-1}}\setminus B_{1}. 
\] 
\end{proof}    
\begin{proof}[Proof of Theorem \ref{thmmain} for $n=6$]   
When $n=6$ and $r=\br=3/2$, we have:  
\begin{align*} 
&\|U\|_{\stpr,\mua^{-1}}\le C(\log\frac{1}{\mua})^{2/3}\\  
&\|U^{r-1}\|_{\stpr,\mua^{-1}}\le C\mua^{-2}\\ 
&\|U^{2^*-2}\|_{\stpr,\mua^{-1}}=\|U\|_{\stpr,\mua^{-1}}\le C(\log\frac{1}{\mua})^{2/3}. 
\end{align*} 
We know
\[  
|Y_{\hg}(\widetilde\xi_{\ya,\lda}^{\hg})-K^{-2}|\le C\mua^4\log\frac1\mua.  
\]    
Hence, (ii) in Corollary \ref{corlowerbound} implies:  
\begin{align}  
\label{6alphaupperbound}  
\alpha\|\ua\|_{L^{\br}(M)}^2\le  
C\Big(\mua^4\big(\log\frac{1}{\mua}\big)^{4/3}  
+\epsa\big(\log\frac{1}{\mua}\big)^{2/3} 
+\mua^6\big(\log\frac{1}{\mua}\big)^{4/3}\Big)  
\end{align}   
From the uniform estimate \eqref{unifest} we derive:    
\begin{align*}  
\epsa\le C\alpha\mua^4\big(\log\frac{1}{\mua}\big)^{1/3},  
\end{align*}  
and by 
Proposition \ref{proplowerbound} we have 
\[ 
\|\ua\|_{L^{\br}(\tBa)}\ge C^{-1}\mua^2\big(\log\frac{1}{\mua}\big)^{2/3}.  
\] 
Inserting   
into \eqref{6alphaupperbound},   
we obtain  
\begin{align*}  
\alpha\mua^4\big(\log\frac{1}{\mua}\big)^{4/3}\le  
C\Big(\mua^4\big(\log\frac{1}{\mua}\big)^{4/3}  
+\alpha\mua^4\log\frac{1}{\mua}\Big).  
\end{align*}  
Once again we obtain $\alpha\le C$, a contradiction.    
Theorem \ref{thmmain} is thus established in the remaining limit case $n=6$.  
\end{proof}    
\section{Appendix: A local to global argument}   
In this Appendix we provide a proof of the local to global argument  
used in Theorem \ref{thmflat}.  
We adapt some ideas from \cite{AL}.   
Let $(M,g)$ be a smooth compact Riemannian manifold
without boundary, $n\ge3$.   
\begin{Lemma}   
\label{lemlocaltoglobal}  
Suppose that there exist $\beps>0$ and $A_{\beps}>0$   
such that   
\begin{align}   
\label{local}   
&\|u\|_{L^{2^*}(M,g)}^2\le K^2\intM\{|\nabg u|^2+c(n)\Rg u^2\}\dvg   
+A_{\beps}\|u\|_{L^1(M,g)}^2,   
\end{align}   
for all $u\in H^1(M)$ such that ${\rm diam}_g({\rm supp})\,u<\beps$.  
Then there exists a constant $A>0$ such that   
\begin{align}   
\label{global}   
&\|u\|_{L^{2^*}(M,g)}^2\le K^2\intM\{|\nabg u|^2+c(n)\Rg u^2\}\dvg   
+A\|u\|_{L^1(M,g)}^2,   
&&\forall u\in H^1(M).   
\end{align}   
\end{Lemma}   
\begin{proof}   
By contradiction.   
Suppose \eqref{global} is not true.  
Then by density of smooth functions in $H^1(M)$,  
for all $\alpha>0$ there exists $\ra\in(1,2)$  
such that  
\begin{align}  
\label{loctoglobcontrassump}  
\la:=\inf_{u\in H^1(M)\setminus\{0\}}  
\frac{\intM\{|\nabg u|^2+c(n)\Rg u^2\}\dvg+\alpha\|u\|_{L^{\ra}(M)}^2}  
{\|u\|_{L^{2^*}(M)}^2}<K^{-2}.  
\end{align}  
By the results in Section \ref{secprelims} with $r=\ra$,  
there exists $\ua\in H^1(M)$, $\ua\ge0$, $\intM\ua^{2^*}\dvg=1$  
such that $\la=\Ia(\ua)$.  
Moreover, $\ua$ satisfies the Euler-Lagrange equation:  
\begin{align}   
\label{appuaeq}   
-\Delg\ua+c(n)\Rg\ua+\alpha\|\ua\|_{L^{\ra}(M)}^{2-\ra}\ua^{\ra-1}=\la\ua^{2^*-1}   
&&\text{on}\;M.   
\end{align}   
Denote by $\xa$ a maximum point of $\ua$.   
By Corollary \ref{coronepointconc}, $\ua$ concentrates   
in energy at $\xa$.   
In particular, for any fixed $\eps>0$,   
\begin{align*}   
\lim_{\alpha\to+\infty}  
\int_{M\setminus B_{\eps}(\xa)}\{|\nabg\ua|^2+\ua^{2^*}\}\dvg=0.  
\end{align*}   
For a fixed $0<\eps<\beps/9$, denote by $\eta$ a smooth cutoff function    
such that   
$\eta\equiv1$ in $B_{2\eps}(\xa)$, $\eta\equiv0$ in $M\setminus B_{4\eps}(\xa)$,   
$0\le\eta\le1$, $|\nabg\eta|\le\eps^{-1}$ in $M$.   
Then, by \eqref{local} and the H\"older inequality,   
\begin{align*}   
\|\eta\ua\|_{L^{2^*}(M,g)}^2  
\le K^2\intM\{|\nabg(\eta\ua)|^2+&c(n)\Rg(\eta\ua)^2\}\dvg\\  
+&A_{\beps}(\volg M)^{2-2/\ra}\|\eta\ua\|_{L^{\ra}(M)}^2,  
\end{align*}   
and consequently,   
\begin{align*}   
\|\ua\|_{L^{2^*}(B_{2\eps}(\xa))}^2   
\le K^2\intM\{|\nabg\ua|^2+&c(n)\Rg\ua^2\}\dvg   
+C\|\ua\|_{L^{\ra}(M)}^2\\   
\nonumber   
&+C\int_{B_{4\eps}(\xa)\setminus B_{2\eps}(\xa)}\{|\nabg\ua|^2+\ua^2\}\dvg.   
\end{align*}   
In turn, using the contradiction assumption \eqref{loctoglobcontrassump}, we have   
\begin{align*}   
\|\ua\|_{L^{2^*}(B_{2\eps}(\xa))}   
\le&K^2\la-(\alpha K^2-C)\|\ua\|_{L^{\ra}(M)}^2\\  
&+C\int_{B_{4\eps}(\xa)\setminus B_{2\eps}(\xa)}\{|\nabg\ua|^2+\ua^2\}\dvg.   
\end{align*}  
Using the expansion  
\[   
\|\ua\|_{L^{2^*}(B_{2\eps}(\xa))}^2   
=1-O(1)\|\ua\|_{L^{2^*}(M\setminus B_{2\eps}(\xa))}^{2^*},   
\]  
and recalling that $\la K<1$,   
we obtain   
\begin{align*}   
\alpha\|\ua\|_{L^{\ra}(M)}^2\le C\|\ua\|_{L^{2^*}(M\setminus B_{2\eps}(\xa))}^2   
+C\int_{B_{4\eps}(\xa)\setminus B_{2\eps}(\xa)}\{|\nabg\ua|^2+\ua^2\}\dvg.   
\end{align*}   
Now let $\eta$ be a cutoff function supported in $M\setminus B_\eps(\xa)$.   
Multiplying \eqref{appuaeq} by $\eta^2\ua$ and integrating by parts,   
we find   
\[   
\intM\eta^2|\nabg\ua|^2\dvg   
\le C\int_{\text{supp}\eta}(\ua^2+\ua^{2^*})\dvg.   
\]   
Therefore,   
\begin{align}   
\label{localineq}   
\alpha\|\ua\|_{L^{\ra}(M)}^2\le C\big(\|\ua\|_{L^2(M\setminus B_{\eps}(\xa))}^2   
+\|\ua\|_{L^{2^*}(M\setminus B_{\eps}(\xa))}^2\big).   
\end{align}   
Finally, by Moser iterations,  
\begin{align}   
\label{moserest}   
\|\ua\|_{L^\infty(M\setminus B_\eps)}   
\le C\|\ua\|_{L^1(M)}\le C(\volg M)^{1-1/\ra}\|\ua\|_{L^{\ra}(M)},   
\end{align}   
see Corollary \ref{coreasysupest}.  
The estimates \eqref{localineq}--\eqref{moserest}   
imply $\alpha\le C$, a contradiction,  
and \eqref{global} is established.  
\end{proof}   
\paragraph{Acknowledgment}  
The second author is grateful to the Rutgers University Mathematics Department  
for hospitality, during her   
leave of absence from Universit\`a di Napoli Federico II.  
    
\end{document}